\documentclass[a4paper,11pt]{article}
\usepackage{amsmath,amssymb,amsthm}
\usepackage{enumitem}
\usepackage{hyperref,url}
\usepackage{color}
\usepackage[margin=30mm]{geometry}
\usepackage{graphicx,subfigure}

\newcommand{\N}{\mathbb{N}}

\newcommand{\II}{\mathbb{I}}
\newcommand{\R}{\mathbb{R}}
\newcommand{\KK}{\mathcal{K}}
\newcommand{\NN}{\mathcal{N}}
\newcommand{\dx}{\, dx}
\newcommand{\dt}{\, dt}
\newcommand{\dd}{\, d}
\newcommand{\pd}{\partial}
\newcommand{\pdnu}{\pd_{n}}
\newcommand{\abs}[1]{\left| #1 \right|}
\newcommand{\norm}[1]{\| #1 \|}
\newcommand{\inner}[2]{\langle #1 , #2 \rangle}
\newcommand{\eps}{\varepsilon}
\newcommand{\Laplace}{\Delta}
\newcommand{\mean}[1]{\overline{ #1 }}

\theoremstyle{plain}
\newtheorem{thm}{Theorem}[]

\newtheorem{remark}{Remark}[section]

\numberwithin{equation}{section}

\title{Cahn--Hilliard inpainting with the double obstacle potential}

\author{Harald Garcke \footnotemark[1] \and Kei Fong Lam \footnotemark[2] \and Vanessa Styles \footnotemark[3]}

\date{ }

\begin{document}

\maketitle

\renewcommand{\thefootnote}{\fnsymbol{footnote}}
\footnotetext[1]{Fakult\"at f\"ur Mathematik, Universit\"at Regensburg, 93040 Regensburg, Germany
({\tt Harald.Garcke@mathematik.uni-regensburg.de}).}
\footnotetext[2]{Department of Mathematics, The Chinese University of Hong Kong, Shatin, N.T., Hong Kong ({\tt kflam@math.cuhk.edu.hk}).}
\footnotetext[3]{University of Sussex, Sussex House, Falmer, Brighton, BN1 9RH, United Kingdom
({\tt V.Styles@sussex.ac.uk}).}

\begin{abstract}
The inpainting of damaged images has a wide range of applications, and many different mathematical methods have been proposed to solve this problem.  Inpainting with the help of Cahn--Hilliard models has been particularly successful, and it turns out that Cahn--Hilliard inpainting with the double obstacle potential can lead to better results compared to inpainting with a smooth double well potential.  However, a mathematical analysis of this approach is missing so far.  In this paper we give first analytical results for a Cahn--Hilliard double obstacle inpainting model regarding existence of global solutions to the time-dependent problem and stationary solutions to the time-independent problem without constraints on
the parameters involved.  With the help of numerical results we show the effectiveness
of the approach for binary and grayscale images.
\end{abstract}

\noindent \textbf{Key words. } Inpainting, Cahn--Hilliard model, double obstacle potential, binary and grayvalue images \\

\noindent \textbf{AMS subject classification. } 49J40, 94A08, 68U10, 35K55

\section{Introduction}

The action of restoring details in missing or damaged portions of an image, commonly known as inpainting, is an active area of mathematical image processing for which several celebrated algorithms have been proposed.  Analogous to disocclusion, i.e., the recovery of scenic information obstructed by visible objects, it is desirable to produce a result which, to the naked eye, does not distinguish where the inpainting has taken place, and a minimum requirement is to have continuity of both the image intensity and the directions of isophotes (lines of equal grayvalue) across the boundary of the inpainting region.  A first approach, proposed in \cite{Bertalmio2}, employs a third-order nonlinear partial differential equation that propagates information (so-called image smoothness measured by the Laplacian of the intensity) into the inpainting region from its surroundings along the direction of least change in grayvalues.  A few iterations of anisotropic diffusion are recommended after a few iterations of inpainting to retain sharpness of edges in the inpainting.  Further analogy to incompressible Navier--Stokes flow is made clear in the follow-up paper of \cite{Bertalmio1} and offers an explanation for the necessity of diffusion.

Earlier algorithms developed for disocclusion have their origin in a variational framework \cite{MM,NMS}.  Building on the success of variational methods for image segmentation (based on the Mumford--Shah functional \cite{MS}) and image restoration (based on the total variation functional \cite{ROF}), several authors have suggested performing inpainting with variational methods by including an additional fidelity term, which serves to keep the solution close to the given image outside the inpainting region \cite{CKS,CSCCD,ChanShen,ES}.  A good overview of the subject can be found in \cite{Tauber}.

In this work, we are interested in the phase field approach proposed by Bertozzi, Esedo\={g}lu, and Gillette \cite{BEG} for binary image inpainting.  Let us fix the setting and introduce our notation.  In a bounded domain $\Omega \subset \R^{d}$ we have a binary image $I : \Omega \to \{-1, 1\}$ such that the region $\{I = 1\}$ represents the black pixels and the region $\{I = -1 \}$ represents the white pixels.  It is assumed that the image has been damaged in a subregion $D \subset \Omega$.  In particular, knowledge about the values of $I(x)$ for $x \in D$ is lost.  The image $I$ is not damaged in the complement set $\Omega \setminus D$.  The approach of Bertozzi, Esedo\={g}lu and Gillette \cite{BEG} involves the following modified Cahn--Hilliard equation with fidelity term:
\begin{subequations}\label{CH}
\begin{alignat}{3}
\pd_{t} u  & = \Laplace w + \lambda (I - u) , \label{CH:1} \\
w & = - \eps \Laplace u + \eps^{-1} W'(u), \label{CH:2}
\end{alignat}
\end{subequations}
where $I$ is the image, $W$ is a potential with equal minima at $\pm 1$, $u: \Omega \to [-1,1]$ can be viewed as the image intensity with $\{u = 1\}$ representing black pixels and $\{u = -1 \}$ representing white pixels, and $\eps > 0$ is a small parameter associated to the thickness of the interfacial layers $\{-1 < u < 1 \}$ separating the black and white regions.  The function $\lambda : \Omega \to \R$ is defined as
\begin{align*}
\lambda(x) = \begin{cases} \alpha & \text{ if } x \in \Omega \setminus D, \\
0 & \text{ if } x \in D,
\end{cases}
\end{align*}
where $\alpha$ is a positive constant. For large values of $\alpha$, the function $\lambda$ penalizes large deviations of the recovery image from the original image in the undamaged regions. 

One of the chief advantages of inpainting with the modified Cahn--Hilliard equation \eqref{CH} compared to other PDE-based approaches, as demonstrated numerically in \cite{BEG}, is the considerable speed-up in computation time thanks to the fast numerical techniques available for the Cahn--Hilliard equation.  In terms of the mathematical analysis of the modified Cahn--Hilliard equation \eqref{CH}, first results concerning global existence and uniqueness of two-dimensional weak solutions to the time-dependent problem and asymptotic behaviour of possible stationary solutions in the limit $\alpha \to \infty$ can be found in \cite{BEGanal}.  In particular, in \cite{BEGanal} the authors established the connection between the approach of \cite{Bertalmio2} that prefers to impose boundary conditions at the edge of the inpainting region $D$, and the variational approach of \cite{ChanShen} that prefers to use a fidelity term.  The existence of weak solutions to the stationary problem is later addressed in \cite{BHS} provided the penalization parameter satisfies $\alpha \geq O(\eps^{-3})$.  Asymptotic behaviour in time, in terms of finite-dimensional attractors, has been established in \cite{CFMattractor}.  These properties of the PDE problem are important in providing a mathematically sound foundation for the development of future algorithms.

In the above contributions, the potential $W$ is taken to be the classical smooth quartic double well function $W_{\mathrm{qu}}(s) = (s^{2}-1)^{2}$, which permits the authors to express \eqref{CH} as one equation in the variable $u$:
\begin{align}\label{onequ}
\pd_{t}u = -\eps \Laplace^{2} u + \eps^{-1} \Laplace W_{\mathrm{qu}}'(u) + \lambda (I - u).
\end{align}
It is well known in the phase field community that solutions $u$ to the Cahn--Hilliard equation (\eqref{CH} with $\lambda = 0$) for a smooth double well potential can attain values outside the physically relevant interval $[-1,1]$; see for instance \cite[Remark 2.1]{CMZ}.  For many practical applications, it is not meaningful to have $u > 1$ or $u < -1$.  One remedy is to employ an ad hoc projection at each numerical iteration  to project any values larger than $1$ to $1$ and values smaller than $-1$ to $-1$.  Another way is to use smaller values of the parameter $\eps$ during the implementation, which has the effect of reducing the deviation of $u$ from $[-1,1]$.  However, neither is an attractive option in the view of implementation as the first remedy can cause small perturbations to propagate throughout the implementation, while the second remedy can significantly increase computational effort as the underlying mesh (for example, for a finite element method) needs to be rather fine for small $\eps$ to have enough elements to resolve all the details.  A third remedy is to use nonsmooth potentials of logarithmic type, as originally proposed by Cahn and Hilliard:
\begin{align}\label{logpot}
W_{\mathrm{log}}(u) = \frac{\theta}{2} \left ( (1+u) \log (1+u) + (1-u) \log(1-u) \right ) + \frac{\theta_{c}}{2} (1-u^{2}),
\end{align}
where $0 < \theta < \theta_{c}$.  The singular behaviour of the derivative $W_{\mathrm{log}}'(s)$ at $s = \pm 1$ enforces the bounds $u \in (-1,1)$.  Therefore, numerical methods that preserve this sort of property (for example a finite difference scheme analyzed in \cite{CWWW}), where the discrete solution stays in $(-1,1)$, circumvent any need for ad hoc processing of the recovery image.  For inpainting applications, the Cahn--Hilliard inpainting model \eqref{CH} with a logarithmic potential $W_{\mathrm{log}}$ has been studied in \cite{CFMlog} for the existence of weak solutions, and numerical simulations performed with $W_{\mathrm{log}}$ reach steady states faster than simulations with $W_{\mathrm{qu}}$.  Due to the technical difficulty in showing the spatial mean value of $u$ stays strictly in the open interval $(-1,1)$ for arbitrary reference time intervals, the authors in \cite{CFMlog} can only prove a local-in-time existence result, and so it is not known if the solution eventually blows up after some finite time, nor if the solution converges to an equilibrium (which is the desired recovery image we seek).  Such uncertainties in solution behaviour may not inspire confidence in a numerical implementation of the model.  After this paper has been accepted, the authors were made aware of a recent review paper \cite{Miran} in which the local-in-time result of \cite{CFMlog} has been extended to be global-in-time by a more careful treatment of the spatial mean value of $u$ (see \cite[Remark 4.3]{Miran} for more details).

Another nonsmooth potential to consider is the double obstacle potential proposed by Oono and Puri \cite{OonoPuri} and Blowey and Elliott \cite{BE1} which is defined as
\begin{align}\label{eq:do}
W(s) = \hat{\beta}(s) + \frac{1}{2}(1-s^{2}), \quad \hat{\beta}(s) = \II_{[-1,1]}(s) = \begin{cases}
0 & \text{ if } s \in [-1,1],\\
+\infty & \text{ otherwise},
\end{cases}
\end{align}
where $\II_{A}$ is the indicator function of the set $A$.  The double obstacle potential can be realised as the ``deep quench'' limit $\theta \to 0$ of the logarithmic potential \eqref{logpot} with $\theta_{c}$ fixed.  With this choice, the equation \eqref{CH:2} should be understood as a variational inequality \cite{BE1}.  Let $\KK := \{ f \in H^{1}(\Omega) : \abs{f} \leq 1 \text{ a.e. in } \Omega \}$; then over a reference time interval $(0,T)$ we ask $(u, w)$ to satisfy
\begin{subequations}\label{Vineq}
\begin{alignat}{3}
\pd_{t} u = \Laplace w + \lambda(x)(I - u) & \text{ in } \Omega \times (0,T), \\
- \left ( w + \eps^{-1} u, v - u \right ) + \left ( \eps \nabla u, \nabla (v - u) \right ) \geq 0 & \text{ for all } v \in \KK, \label{Vineq:2} \\
\pdnu w = 0 & \text{ on } \pd \Omega \times (0,T), \\
u(0) = u_{0} & \text{ in } \Omega,
\end{alignat}
\end{subequations}
where $(\cdot, \cdot)$ denotes the $L^{2}$-scalar product and $\pdnu f$ is the normal derivative of $f$ on $\pd \Omega$.  Implicitly in \eqref{Vineq:2} we have also supplemented the boundary condition $\pdnu u = 0$ on $\pd \Omega \times (0,T)$.  Due to the variational inequality \eqref{Vineq:2} it is not feasible to express the inpainting model \eqref{Vineq} as one equation like \eqref{onequ}.  Consequently, it is unknown if the asymptotic behaviour $\alpha \to \infty$ in matching isophotes proved in \cite{BEGanal} can be reproduced for \eqref{Vineq}.  

The idea of using the double obstacle potential for Cahn--Hilliard inpainting is first suggested in \cite{Bosch}, where the focus is fast and practical solvers for the Moreau--Yosida approximation of \eqref{Vineq}; see equation \eqref{beta:approx} and the system \eqref{CHdelta} below.  It is also reported in the same paper that the double obstacle Cahn--Hilliard inpainting leads to better visual results compared to the variant \eqref{onequ} with the smooth potential and other higher-order TV-based inpainting models provided by Sch\"{o}nlieb \cite{BHS,sb}. However, a mathematical analysis of \eqref{Vineq} is missing so far, and hence to lend weight to the numerical simulations in \cite{Bosch}, our first contribution to the field of Cahn--Hilliard inpainting, albeit an analytical one, is to establish the existence of global weak solutions to \eqref{Vineq}.  Furthermore, in our method of proof, we have to extend the techniques of \cite{BE1} in order to derive estimates in the presence of the data fidelty term $\lambda (I - u)$ in \eqref{CH:1}.  While this type of term is not limited to inpainting, we envision that our methods can be broadly applied to other models of similar structures.

Our second contribution to Cahn--Hilliard inpainting is the existence of strong solutions to the stationary problem of \eqref{Vineq} without constraints on $\alpha$ and $\eps$, as opposed to the restriction $\alpha \geq O(\eps^{-3})$ specified in \cite{BHS} for the quartic double well potential $W_{\mathrm{qu}}$.  As the stationary solution is the desired recovery image in the inpainting process, the existence of such a solution shows that the possible steady states attained in numerical implementations are genuine outputs for the inpainting problem.  

A third contribution is the proposal of a finite element scheme to solve the variational inequality \eqref{Vineq} directly, as opposed to the one used by \cite{Bosch} involving the Moreau--Yosida approximation \eqref{CHdelta}.  We defer to Sec.~\ref{sec:num} for a short discussion on the merits of our chosen numerical approach.

Lastly, we note that there have been some attempts of generalizing the modified Cahn--Hilliard equation \eqref{CH} from binary to grayscale images.  One approach is to split the grayscale image bitwise into channels and apply the binary Cahn--Hilliard inpainting to each binary channel \cite{sb} (which is also done here). Another is to treat the grayvalues of image data as the real part of a complex variable and employ the complex-valued Cahn--Hilliard equation \cite{CFMcomplex, GS}.  A third approach, motivated by the use of total variation for grayvalue image decomposition and restoration \cite{LL,OSV}, is the TV-$H^{-1}$ inpainting method proposed by \cite{BHS}.  But a more natural approach would be to employ a multiphase Cahn--Hilliard system with a fidelity term.  This has been done for smooth multiwell potentials in \cite{Boschfrac, CFMVec}.  While in principle the extension for the double obstacle potential can be made, we defer the formulation, mathematical analysis, and numerical implementation of the multiphase case to future research.

This paper is organised as follows:  In Sec.~\ref{sec:timedep} we study the time-dependent problem \eqref{Vineq} and prove the existence of global weak solutions.  The existence of a strong solution to the stationary problem is addressed in Sec.~\ref{sec:stat}, and in Sec.~\ref{sec:num} we outline a finite element scheme different to the one used by \cite{Bosch} and present several inpainting results for binary and grayscale images.

\paragraph{Notation.}  For convenience, we will often use the notation $L^{p} := L^{p}(\Omega)$ and $W^{k,p} := W^{k,p}(\Omega)$ for any $p \in [1,\infty]$, $k > 0$ to denote the standard Lebesgue spaces and Sobolev spaces equipped with the norms $\norm{\cdot}_{L^{p}}$ and $\norm{\cdot}_{W^{k,p}}$.  By $(\cdot,\cdot)$ we denote the $L^{2}$-inner product.  In the case $p = 2$ we use $H^{k} := W^{k,2}$ and the norm $\norm{\cdot}_{H^{k}}$.  We denote by $H^{2}_{n}(\Omega) := \{ f \in H^{2}(\Omega) : \pdnu f = 0 \text{ on } \pd \Omega \}$.  Moreover, the dual space of a Banach space $X$ will be denoted by $X'$, and the duality pairing between $X$ and $X'$ is denoted by $\inner{\cdot}{\cdot}_{X}$.  

\section{The time-dependent problem}\label{sec:timedep}
We make the following assumptions.
\begin{enumerate}[label=$(\mathrm{A \arabic*})$, ref = $\mathrm{A \arabic*}$]
\item \label{ass:Omega} We consider a bounded domain $\Omega \subset \R^{d}$, $(d=1,2,3)$, with $C^{1,1}$ boundary $\pd \Omega$ or a convex bounded domain, and the damaged region $D \subsetneq \Omega$ is measurable.
\item \label{ass:I} The original image $I \in L^{\infty}(\Omega)$ satisfies $\abs{I(x)} \leq 1$ for a.e. $x \in \Omega$, and we assume that $I$ is not identically equal to $1$ or $-1$ in $\Omega \setminus D$; i.e., we exclude the cases $I \equiv 1$ and $I \equiv -1$ a.e. in $\Omega \setminus D$.
\item \label{ass:u0} $u_{0} \in H^{1}(\Omega)$ and $u_{0}(x) \in [-1,1]$ for a.e. $x \in \Omega$.  
\end{enumerate}

\begin{remark}
In practice a binary image $I$ would take values only in $\{ \pm 1 \}$.  But for the coming analysis, it suffices to assume that $\abs{I(x)} \leq 1$ a.e. in $\Omega$.  Furthermore, the condition that $I$ is not identically equal to $1$ or $-1$ in the undamaged region implies that the undamaged part is not completely black or white which would lead to the trivial recovery of $u \equiv 1$ or $u \equiv - 1$.
\end{remark}

Our result concerning the existence of a solution to the time-dependent problem \eqref{Vineq} is formulated as follows.

\begin{thm}\label{thm:1}
Let $T$, $\alpha$, and $\eps$ be arbitrary positive constants.  Under assumptions \eqref{ass:Omega}, \eqref{ass:I}, and \eqref{ass:u0}, there exists at least one pair $(u, w)$ such that 
\begin{align*}
u \in L^{\infty}(0,T;H^{1}(\Omega)) \cap L^{2}(0,T;H^{2}_{n}(\Omega)) \cap H^{1}(0,T;(H^{1}(\Omega))'), \quad w \in L^{2}(0,T;H^{1}(\Omega)),
\end{align*}
satisfying $u(0) = u_{0}$ a.e. in $\Omega$, $\abs{u} \leq 1$ a.e. in $\Omega \times (0,T)$, and 
\begin{subequations}\label{weakform}
\begin{alignat}{3}
\inner{\pd_{t}u(t)}{\zeta}_{H^{1}} + ( \nabla w(t), \nabla \zeta)  -( \lambda (I - u(t)),  \zeta ) & = 0, \label{weak:1} \\
- \left ( w(t) + \eps^{-1} u(t), v - u(t) \right ) + \left ( \eps \nabla u(t), \nabla (v - u(t)) \right ) & \geq 0 \label{weak:2}
\end{alignat}
\end{subequations}
for all $v \in \KK$, $\zeta \in H^{1}(\Omega)$ and for a.e. $t \in (0,T)$.
\end{thm}

\begin{remark}[Uniqueness]
For Dirichlet boundary conditions $u = -1$ and $w = 0$ on $\pd \Omega \times (0,T)$, uniqueness of solutions is shown in \cite[\S 3.2.1]{Poole}.  The difficulty with the Neumann boundary condition is that the spatial mean of the difference of two solutions may not be zero.  
\end{remark}

The proof will be based on an appropriate approximation of the subdifferential $\pd \II_{[-1,1]}(s)$ and then deriving uniform a priori estimates.  For this purpose, we introduce the classical approximation (see also \cite{Bosch}),
\begin{align}\label{beta:approx}
\beta_{\delta}(s) := \frac{1}{\delta} \left ( \max (0, s-1) + \min(0, s+1) \right ) \text{ for } \delta > 0, \, s \in \R,
\end{align}
and set $\hat{\beta}_{\delta}$ to be the antiderivative of $\beta_{\delta}$ with $\hat{\beta}_{\delta}(0) = 0$; i.e., $\hat{\beta}_{\delta}'(s) = \beta_{\delta}(s)$ and integrating \eqref{beta:approx} leads to the expression
\begin{align}
\label{hat:beta:approx}
\hat{\beta}_{\delta}(s) = \begin{cases}
\frac{1}{2 \delta} (s-1)^{2} & \text{ if } s \geq 1, \\
0 & \text{ if } \abs{s} \leq 1, \\
\frac{1}{2 \delta} (s+1)^{2} & \text{ if } s \leq -1.
\end{cases}
\end{align}
For each $\delta$, we have that $\hat{\beta}_{\delta}$ is nonnegative, convex with quadratic growth, and $\beta_{\delta}$ satisfies
\begin{align}\label{betadelta}
\beta_{\delta}(s) = 0 \text{ for } \abs{s} \leq 1, \quad s \, \beta_{\delta}(s) \geq 0 \text{ for } \abs{s} \geq 1.
\end{align}
We also define the regularized potential
\begin{align}\label{reg:energy}
W_{\delta}(s) = \hat{\beta}_{\delta}(s) + \frac{1}{2}(1-s^{2}).
\end{align}

\subsection{Approximation scheme}
We consider a sequence of solutions $\{(u_{\delta}, w_{\delta})\}_{\delta > 0}$, indexed by $\delta > 0$, such that
\begin{subequations}\label{CHdelta}
\begin{alignat}{3}
\pd_{t} u_{\delta} = \Laplace w_{\delta} + \lambda(x)(I - u_{\delta}) & \text{ a.e. in } \Omega \times (0,T),  \label{CHd:1}\\
w_{\delta}  = - \eps \Laplace u_{\delta} + \eps^{-1} \beta_{\delta}(u_{\delta}) - \eps^{-1} u_{\delta} \, & \text{ a.e. in } \Omega \times (0,T),  \label{CHd:2} \\
\pdnu u_{\delta}  = \pdnu w_{\delta} = 0 & \text{ on } \pd \Omega \times (0,T), \\
u_{\delta}(0)  = u_{0} & \text{ in } \Omega.
\end{alignat}
\end{subequations}
For each $\delta > 0$, the existence and uniqueness of $(u_{\delta}, w_{\delta})$ to \eqref{CHdelta} with the regularities
\begin{align*}
u_{\delta} \in L^{\infty}(0,T;H^{2}(\Omega)) \cap H^{1}(0,T;L^{2}(\Omega)), \quad w_{\delta} \in L^{\infty}(0,T;L^{2}(\Omega)) \cap L^{2}(0,T;H^{2}(\Omega))
\end{align*}
can be proved using standard methods.  We postpone the existence proof of $(u_{\delta},w_{\delta})$ to Sec.~\ref{sec:Approx:Exist}.  The goal is to derive a priori estimates that are uniform in $\delta$ and then pass to the limit $\delta \to 0$.  In the following, the symbol $C$ will denote positive constants independent of $\delta$ and may vary from line to line.

\subsection{Uniform estimates}
\paragraph{First estimate.}
Testing \eqref{CHd:1} with $u_{\delta}$, testing \eqref{CHd:2} with $\Laplace u_{\delta}$, and summing the resulting equations yields
\begin{align*}
& \frac{\dd}{\dt} \frac{1}{2} \norm{u_{\delta}}_{L^{2}}^{2} + \int_{\Omega} \frac{1}{\eps} \underbrace{ \beta_{\delta}'(u_{\delta}) \abs{\nabla u_{\delta}}^{2} }_{\geq 0} + \, \eps \abs{\Laplace u_{\delta}}^{2} \dx  = \int_{\Omega} \lambda (I - u_{\delta}) u_{\delta} + \frac{1}{\eps} \abs{\nabla u_{\delta}}^{2} \dx.
\end{align*}
Thanks to the $\pdnu u_{\delta} = 0$ on $\pd \Omega \times (0,T)$, integrating by parts leads to the inequality
\begin{align}\label{inequDelta}
\norm{\nabla u_{\delta}}_{L^{2}}^{2} \leq \eta \norm{\Laplace u_{\delta}}_{L^{2}}^{2} + \frac{1}{4 \eta} \norm{u_{\delta}}_{L^{2}}^{2} \quad \text{ for all } \eta > 0,
\end{align}
and Gronwall's inequality gives
\begin{align}\label{Apri:1}
\norm{u_{\delta}}_{L^{\infty}(0,T;L^{2})} + \norm{\nabla u_{\delta}}_{L^{2}(0,T;L^{2})} + \norm{\Laplace u_{\delta}}_{L^{2}(0,T;L^{2})} \leq C.
\end{align}

\paragraph{Second estimate.}
Testing \eqref{CHd:1} with $- \eps \Laplace u_{\delta}$ and also with $\frac{1}{\eps} \beta_{\delta}(u_{\delta})$ and testing \eqref{CHd:2} with $- \Laplace w_{\delta}$ give the following identities:
\begin{equation*}
\begin{alignedat}{3}
\frac{\dd}{\dt} \frac{\eps}{2} \norm{\nabla u_{\delta}}_{L^{2}}^{2} &= \int_{\Omega} - \eps \Laplace w_{\delta} \Laplace u_{\delta} - \eps \lambda (I - u_{\delta}) \Laplace u_{\delta} \dx, \\
\frac{\dd}{\dt} \int_{\Omega} \frac{1}{\eps} \hat{\beta}_{\delta}(u_{\delta}) \dx & = - \int_{\Omega} \frac{1}{\eps} \beta_{\delta}'(u_{\delta}) \nabla w_{\delta} \cdot \nabla u_{\delta} - \frac{1}{\eps} \lambda (I - u_{\delta}) \beta_{\delta}(u_{\delta}) \dx, \\
\norm{\nabla w_{\delta}}_{L^{2}}^{2} & = \int_{\Omega} \frac{1}{\eps} \beta_{\delta}'(u_{\delta}) \nabla u_{\delta} \cdot \nabla w_{\delta} - \frac{1}{\eps} \nabla u_{\delta} \cdot \nabla w_{\delta} + \eps \Laplace u_{\delta} \Laplace w_{\delta} \dx.
\end{alignedat}
\end{equation*}
Upon summing we obtain
\begin{align*}
& \frac{\dd}{\dt} \left ( \frac{\eps}{2} \norm{\nabla u_{\delta}}_{L^{2}}^{2} + \frac{1}{\eps} \norm{\hat{\beta}_{\delta}(u_{\delta})}_{L^{1}} \right ) + \norm{\nabla w_{\delta}}_{L^{2}}^{2} \\
& \quad = \int_{\Omega} - \frac{1}{\eps} \nabla w_{\delta} \cdot \nabla u_{\delta} - \eps \lambda (I - u_{\delta}) \Laplace u_{\delta} + \frac{1}{\eps} \lambda (I - u_{\delta}) \beta_{\delta}(u_{\delta})  \dx \\
& \quad =: J_{1} + J_{2} + J_{3}.
\end{align*}
The first and second terms on the right-hand side can be estimated using Young's inequality
\begin{align*}
J_{1} + J_{2} \leq \frac{1}{2} \norm{\nabla w_{\delta}}_{L^{2}}^{2} + \frac{1}{2 \eps^{2}} \norm{\nabla u_{\delta}}_{L^{2}} + \eps \alpha \left ( 1 + \norm{u_{\delta}}_{L^{2}} \right ) \norm{\Laplace u_{\delta}}_{L^{2}}.
\end{align*}
For the third term, we claim that
\begin{align}\label{I-u:beta}
(I - u_{\delta}) \beta_{\delta}(u_{\delta})  \begin{cases}
= 0 & \text{ if } \abs{u_{\delta}} \leq 1, \\
\leq 0 & \text{ if } \abs{u_{\delta}} > 1
\end{cases}
\end{align}
follows from \eqref{betadelta}.  Indeed, for the first case we use $\beta_{\delta}(s) = 0$ for $\abs{s} \leq 1$, and for the second case, suppose at some point $x \in \Omega$, $u_{\delta}(x) > 1$.  Then, $\beta_{\delta}(u_{\delta}(x)) \geq 0$ and $(I(x) - u_{\delta}(x)) \leq 0$, which yields that the product is nonpositive.  A similar argument also applies for the case $u_{\delta}(x) < -1$.  Then, the integral $J_{3}$ is nonpositive, which leads to
\begin{align*}
\frac{\dd}{\dt} \left ( \norm{\nabla u_{\delta}}_{L^{2}}^{2} + \norm{\hat{\beta}_{\delta}(u_{\delta})}_{L^{1}} \right ) + \norm{\nabla w_{\delta}}_{L^{2}}^{2} \leq C \left ( 1 + \norm{u_{\delta}}_{H^{1}}^{2}  + \norm{\Laplace u_{\delta}}_{L^{2}}^{2} \right ),
\end{align*}
and by Gronwall's inequality and \eqref{Apri:1} we obtain
\begin{align}\label{Apri:2}
\norm{\hat{\beta}_{\delta}(u_{\delta})}_{L^{\infty}(0,T;L^{1})} + \norm{\nabla u_{\delta}}_{L^{\infty}(0,T;L^{2})} + \norm{\nabla w_{\delta}}_{L^{2}(0,T;L^{2})} \leq C.
\end{align}
Let us mention that since $u_{\delta}(0) = u_{0} \in [-1,1]$ a.e. in $\Omega$, by \eqref{hat:beta:approx} it is easy to see that $\hat{\beta}_{\delta}(u_{\delta}(0)) = 0$ a.e. in $\Omega$.

\paragraph{Third estimate.}
By inspection of \eqref{CHd:1} we easily infer the estimate on the time derivative
\begin{align}\label{Apri:3}
\norm{\pd_{t} u_{\delta}}_{L^{2}(0,T;(H^{1})')} \leq \norm{\nabla w_{\delta}}_{L^{2}(0,T;L^{2})} + \alpha \left ( 1 + \norm{u_{\delta}}_{L^{2}(0,T;L^{2})} \right ) \leq C.
\end{align}
Furthermore, for the mean value $\mean{u_{\delta}} := \frac{1}{\abs{\Omega}} \int_{\Omega} u_{\delta} \dx$, we also obtain that
\begin{align*}
\abs{\pd_{t} \mean{u_{\delta}}} = \frac{1}{\abs{\Omega}} \abs{\int_{\Omega \setminus D} \alpha (I - u_{\delta}) \dx} \leq C \left ( 1  + \norm{u_{\delta}}_{L^{2}} \right ) \in L^{\infty}(0,T),
\end{align*}
and so
\begin{align}\label{Apri:4}
\norm{\mean{u_{\delta}}}_{W^{1,\infty}(0,T)} \leq C.
\end{align}

\paragraph{Fourth estimate.}
With the Aubin--Lions lemma we obtain from \eqref{Apri:1}, \eqref{Apri:2}, and \eqref{Apri:3} that there exists a function $u$ such that along subsequences $\{\delta_{k}\}_{k \in \N}$, $\delta_{k} \to 0$ as $k \to \infty$,
\begin{align*}
u_{\delta_{k}} \to u \text{ strongly in }C^{0}([0,T];L^{r}) \text{ and a.e. in } \Omega \times (0,T)
\end{align*}
for any $r < \infty$ for spatial dimensions $d = 1,2$ and for $r < 6$ in three spatial dimensions.  Furthermore, we infer from \eqref{Apri:4} the equiboundedness and equicontinuity of $\{\overline{u_{\delta_{k}}} \}_{k \in \N}$ thanks to the fundamental theorem of calculus:
\begin{align*}
\abs{\overline{u_{\delta_{k}}}(r) - \overline{u_{\delta_{k}}}(s)} = \abs{\int_{s}^{r} \pd_{t} \overline{u_{\delta_{k}}}(t) \dt} \leq C \abs{r-s}.
\end{align*}
By virtue of the Ascoli--Arzel\`a theorem, we infer
\begin{align}\label{str:mean}
\mean{u_{\delta_{k}}}(t) \to \mean{u}(t) \text{ strongly in } C^{0}([0,T]).
\end{align}
We now claim that $u(t) \in [-1,1]$ a.e. in $\Omega$ for all $t \in (0,T)$.  From the definition of $\hat{\beta}_{\delta}$ in \eqref{hat:beta:approx} and the a priori estimate \eqref{Apri:2} we infer that
\begin{align*}
\int_{\{ u_{\delta_{k}} > 1\}} (u_{\delta_{k}}(t) - 1)^{2} \dx + \int_{\{ u_{\delta_{k}} < -1\}} (u_{\delta_{k}}(t) + 1)^{2} \dx \leq 2 \delta_k \norm{\hat{\beta}_{\delta_{k}}(u_{\delta_{k}}(t))}_{L^{1}} \leq C \delta_{k}
\end{align*}
for all $t \in (0,T)$, which is equivalent to
\begin{align*}
\int_{\Omega} (u_{\delta_{k}}(t)-1)_{+}^{2} + (-u_{\delta_{k}}(t) -  1)_{-}^{2} \dx \leq C \delta_{k},
\end{align*}
where $f_{+} := \max (f, 0)$ and $f_{-} := \max (-f, 0)$.  Hence, the limit function $u$ satisfies
\begin{align*}
\int_{\Omega} (-u(t) -1)_{-}^{2} + (u(t) -1)_{+}^{2} \dx = 0 \text{ for all } t \in (0,T),
\end{align*}
which implies that $ \abs{u(t)} \leq 1$ for all $t \in (0,T)$. 

Let $\eta(t) \in C^{\infty}_{c}(0,T)$ be an arbitrary test function.  Then, passing to the limit $\delta_{k} \to 0$ in the equation
\begin{align*}
\int_{0}^{s} \eta(t) \int_{\Omega} \pd_{t} u_{\delta_{k}} \dx \dt = \int_{0}^{s} \eta(t) \int_{\Omega \setminus D} \alpha (I - u_{\delta_{k}}) \dx \dt
\end{align*}
leads to
\begin{align*}
\inner{\pd_{t} u(t)}{1}_{H^{1}} = \int_{\Omega \setminus D} \alpha (I - u(t)) \dx \text{ for all } t \in (0,T).
\end{align*}
We now claim that $\mean{u}(t) \in (-1,1)$ for all $t \in (0,T)$.  Suppose to the contrary there exists a time $t_{*} \in (0,T)$ such that $\mean{u}(t_{*}) = 1$.  This implies that $u(t_{*},x) \equiv 1$ a.e. in $\Omega$, and so
\begin{align*}
\inner{\pd_{t} u(t_{*})}{1}_{H^{1}} = \int_{\Omega \setminus D} \alpha (I - u(t_{*})) \dx < 0
\end{align*}
due to the assumption that $I$ is not identically equal to $1$ or $-1$ in $\Omega \setminus D$.  Hence, the mean $\mean{u}(t)$ must be strictly decreasing in a neighbourhood around $t_{*}$, i.e., $\mean{u}(t) > 1$ for $t < t_{*}$, and this contradicts the fact that $\abs{u(s)} \leq 1$ for all $s \in (0,T)$.  By a similar argument, if $t_{*} \in (0,T)$ is such that $\mean{u}(t_{*}) = -1$, then $u(t_{*},x) \equiv -1$ a.e. in $\Omega$, and as $I - u(t_{*}) > 0$ a.e. in $\Omega \setminus D$ this yields that $\mean{u}(t)$ is strictly increasing in a neighbourhood of $t_{*}$.  Hence $\mean{u}(t) < -1$ for $t < t_{*}$, and this violates the fact that $\abs{u(s)} \leq 1$ for all $s \in (0,T)$.  

By the above reasoning, we find that $\mean{u}(t) \in (-1,1)$ for all $t \in (0,T)$.  We now derive additional a priori estimates for the mean $\mean{w}_{\delta}$ in order to pass to the limit in \eqref{CHd:2}.  Testing \eqref{CHd:2} with $\pm 1$ yields
\begin{align*}
\abs{ \int_{\Omega} w_{\delta} (t) \dx } \leq \int_{\Omega} \frac{1}{\eps} \abs{\beta_{\delta}(u_{\delta}(t))} + \frac{1}{\eps} \abs{u_{\delta}(t)} \dx.
\end{align*}
Then, testing \eqref{CHd:2} with $u_{\delta}$ leads to
\begin{align*}
\eps \norm{\nabla u_{\delta}(t)}_{L^{2}}^{2} + \int_{\Omega} \frac{1}{\eps} u_{\delta}(t) \beta_{\delta}(u_{\delta}(t)) \dx = \int_{\Omega} w_{\delta}(t) u_{\delta}(t) + \frac{1}{\eps} \abs{u_{\delta}(t)}^{2} \dx.
\end{align*}
Combining the above two estimates and using the fact that $\abs{\beta_{\delta}(r)} \leq r \beta_{\delta}(r)$ for all $r \in \R$ we now arrive at
\begin{align}\label{mean:1}
\abs{ \int_{\Omega} w_{\delta} (t) \dx } \leq \frac{C}{\eps} \norm{u_{\delta}(t)}_{L^{2}}^{2} + \int_{\Omega} w_{\delta}(t) u_{\delta}(t)  \dx.
\end{align}
Next, consider the function $f$ that solves $-\Laplace f = (u_{\delta} - \mean{u_{\delta}})(t)$ with homogeneous Neumann boundary conditions and the compatibility condition $\mean{f} = 0$.  By the Lax--Milgram theorem and the Poincar\'{e} inequality we obtain that $f \in H^{1}(\Omega)$ with
\begin{align*}
\norm{f}_{H^{1}} \leq c\norm{u_{\delta}(t) - \mean{u_{\delta}}(t)}_{L^{2}} \leq c^{2} \norm{\nabla u_{\delta}(t)}_{L^{2}},
\end{align*}
where $c > 0$ denotes the constant from the Poincar\'{e} inequality.  Furthermore, testing with $w_{\delta}$ in the variational formulation for $f$ and testing \eqref{CHd:1} with $f$ yields
\begin{equation}\label{mean:2}
\begin{aligned}
 \int_{\Omega} \lambda (I - u_{\delta}(t)) f - \pd_{t} u_{\delta}(t) f \dx & = \int_{\Omega} \nabla f \cdot \nabla w_{\delta}(t) \dx = \int_{\Omega} (u_{\delta}(t) - \mean{u_{\delta}}(t)) w_{\delta}(t) \dx \\
 & = \int_{\Omega} u_{\delta}(t) w_{\delta}(t) - (\mean{u_{\delta}}(t) - \mean{u}(t) + \mean{u}(t)) w_{\delta} \dx.
\end{aligned}
\end{equation}
Substituting \eqref{mean:2} into \eqref{mean:1} and rearranging gives
\begin{equation}\label{mean:3}
\begin{aligned}
& \left ( 1 - \abs{\mean{u}(t)} - \sup_{t \in (0,T)}\abs{\mean{u_{\delta}}(t) - \mean{u}(t)} \right )\abs{\int_{\Omega} w_{\delta}(t) \dx} & \\
& \quad \leq \frac{C}{\eps} \norm{u_{\delta}(t)}_{L^{2}}^{2} + \int_{\Omega} \lambda (I - u_{\delta}(t)) f - \pd_{t} u_{\delta}(t) f \dx \\
& \quad \leq \frac{C}{\eps} \norm{u_{\delta}(t)}_{L^{2}}^{2}  + C \left ( 1 + \norm{u_{\delta}(t)}_{L^{2}} + \norm{\pd_{t}u_{\delta}(t)}_{(H^{1})'} \right ) \norm{\nabla u_{\delta}(t)}_{L^{2}}.
\end{aligned}
\end{equation}
For the subsequence $\{\delta_{k}\}_{k \in \N}$ where \eqref{str:mean} holds, the strong convergence \eqref{str:mean} implies that there exists some $k_{0} > 0$ such that for all $k > k_{0}$
\begin{align*}
\sup_{t \in (0,T)} \abs{\mean{u_{\delta_{k}}}(t) - \mean{u}(t)} \leq \frac{1}{2} \sup_{t \in (0,T)} \left ( 1 - \abs{\mean{u}(t)} \right ).
\end{align*}
Since $\abs{\mean{u}(t)} < 1$ for all $t \in (0,T)$ and since $\mean{u}$ is continuous, the prefactor on the left-hand side of \eqref{mean:3} is bounded away from $0$ uniformly in $t$.  This yields that (for the subsequence $\{\delta_{k}\}_{k \in \N}$)
\begin{align*}
\abs{\mean{w_{\delta_{k}}}} \in L^{2}(0,T),
\end{align*}
and by the Poincar\'{e} inequality 
\begin{align}\label{Apri:5}
\norm{w_{\delta_{k}}}_{L^{2}(0,T;L^{2})} \leq C.
\end{align}
With \eqref{Apri:5} in hand, testing \eqref{CHd:2} with $\beta_{\delta_{k}}(u_{\delta_{k}})$ leads to
\begin{align*}
\frac{1}{\eps} \norm{\beta_{\delta_{k}}(u_{\delta_{k}})}_{L^{2}}^{2} + \int_{\Omega}  \underbrace{\beta_{\delta_{k}}'(u_{\delta_{k}}) \abs{\nabla u_{\delta_{k}}}^{2}}_{\geq 0} \dx \leq \left ( \norm{w_{\delta_{k}}}_{L^{2}} + \frac{1}{\eps} \norm{u_{\delta_{k}}}_{L^{2}} \right ) \norm{\beta_{\delta_{k}}(u_{\delta_{k}})}_{L^{2}},
\end{align*}
which implies that 
\begin{align}\label{Apri:6}
\norm{\beta_{\delta_{k}}(u_{\delta_{k}})}_{L^{2}(0,T;L^{2})} + \norm{u_{\delta_{k}}}_{L^{2}(0,T;H^{2})} \leq C,
\end{align}
where the latter estimate comes from elliptic regularity theory (\cite[Thm.~2.4.2.7]{Grisvard} for bounded domains with $C^{1,1}$-boundary or \cite[Thm.~3.2.1.3]{Grisvard} for general bounded convex domains) applied to \eqref{CHd:2}.

\subsection{Passing to the limit}
From \eqref{Apri:1}, \eqref{Apri:2}, \eqref{Apri:3}, \eqref{Apri:4}, \eqref{Apri:5}, and \eqref{Apri:6} we obtain (for a nonrelabelled subsequence)
\begin{align*}
u_{\delta} & \to u \text{ weakly-* in } L^{\infty}(0,T;H^{1}(\Omega)), \\
u_{\delta} & \to u  \text{ weakly in }  L^{2}(0,T;H^{2}(\Omega)) \cap H^{1}(0,T;(H^{1}(\Omega))'), \\
u_{\delta} & \to u  \text{ strongly in } C^{0}([0,T];L^{r}(\Omega)) \cap L^{2}(0,T;W^{1,r}(\Omega)), \\
w_{\delta} & \to w \text{ weakly in } L^{2}(0,T;H^{1}(\Omega)), \\
\beta_{\delta}(u_{\delta}) & \to \xi  \text{ weakly in } L^{2}(0,T;L^{2}(\Omega))
\end{align*}
to some limit functions $(u, w, \xi)$, where $r < \infty$ for $d = 1,2$ and $r < 6$ for $d = 3$.  Standard results in maximal monotone operator theory lead to the assertion $\xi \in \pd \II_{[-1,1]}(u)$ a.e.~in $\Omega \times (0,T)$.  Testing \eqref{CHd:1} and \eqref{CHd:2} with $\eta(t) \zeta(x)$, where $\eta \in C^{\infty}_{c}(0,T)$ and $\zeta \in H^{1}(\Omega)$ are arbitrary, then passing to the limit yields that $(u,w,\xi)$ satisfies
\begin{subequations}
\begin{alignat}{3}
0 & = \inner{\pd_{t} u(t)}{\zeta}_{H^{1}} + \int_{\Omega} \nabla w(t) \cdot \nabla \zeta - \lambda(x) (I - u(t)) \zeta \dx, \\
0 & = \int_{\Omega} w(t) \zeta - \frac{1}{\eps} \xi(t) \zeta + \frac{1}{\eps} u(t) \zeta - \eps \nabla u(t) \cdot \nabla \zeta \dx, \label{ptl:2}
\end{alignat}
\end{subequations}
for all $\zeta \in H^{1}(\Omega)$ and for a.e. $t \in (0,T)$.  For an arbitrary $v \in \KK$ we have $\II_{[-1,1]}(v) = 0$, and from the definition of $\xi$ belonging to the subdifferential $\pd \II_{[-1,1]}(u)$,
\begin{align*}
0 = \II_{[-1,1]}(v) - \II_{[-1,1]}(u) \geq (\xi, v - u),
\end{align*}
we obtain by substituting $\zeta = v - u(t)$ in \eqref{ptl:2} the variational inequality \eqref{Vineq:2}.  This concludes the proof of Theorem~\ref{thm:1}, and it only remains to show the existence of a solution to the approximation system \eqref{CHdelta}.

\subsection{Existence to approximation problem}\label{sec:Approx:Exist}
The approximation system \eqref{CHdelta} can be seen as a Cahn--Hilliard equation with source terms and a regular potential with quadratic growth.  Indeed, from \eqref{hat:beta:approx} and \eqref{reg:energy}, we see that for fixed $\delta > 0$, both $\hat{\beta}_{\delta}$ and $W_{\delta}$ are $C^{2}$-functions with bounded second derivatives.  Below we will briefly sketch the derivation of the a priori estimates necessary to justify the testing procedures above.  For simplicity, we set $\eps = 1$ and drop the subscript $\delta$.  Employing a Galerkin approximation and the usual testing procedure for the Cahn--Hilliard equation (i.e., combining the equalities obtained from testing \eqref{CHd:1} with $w$ and testing \eqref{CHd:2} with $\pd_{t}u$), which is also possible on the Galerkin level by choosing eigenfunctions of the Laplace operator with Neumann boundary conditions, firstly we obtain
\begin{align}\label{delta:Apri}
\frac{\dd}{\dt} \int_{\Omega} W(u) + \frac{1}{2} \abs{\nabla u}^{2} \dx + \int_{\Omega} \abs{\nabla w}^{2} \dx = \int_{\Omega} \lambda (I - u) w \dx,
\end{align}
where $W$ is the regularized energy as defined in \eqref{reg:energy}.  Since $W'$ has linear growth, a structural assumption of the form
\begin{align*}
\abs{W'(s)}^{2} \leq c_{0} \left ( 1 + W(s) \right ), \quad W(s) \geq c_{1} \abs{s}^{2} - c_{2} \quad \text{ for all } s \in \R, 
\end{align*}
for positive constants $c_{0}, c_{1}, c_{2}$ allows us to estimate the right-hand side of \eqref{delta:Apri} as follows:
\begin{align*}
\abs{\mathrm{RHS}} & \leq C \left ( 1 + \norm{u}_{L^{2}} \right ) \left ( \norm{w - \mean{w}}_{L^{2}} + \abs{\mean{w}} \right ) \leq C \left ( 1 + \norm{u}_{L^{2}}^{2} + \norm{W'(u)}_{L^{1}}^{2} \right ) + \frac{1}{2} \norm{\nabla w}_{L^{2}}^{2} \\
& \leq C \left ( 1 + \norm{W(u)}_{L^{1}} \right ) + \frac{1}{2} \norm{\nabla w}_{L^{2}}^{2},
\end{align*}
where we have used that $\mean{w} = \mean{W'(u)}$ and the Poincar\'{e} inequality.  Estimating the right-hand side of \eqref{delta:Apri} with this inequality and using a Gronwall argument, and then using the estimates for the mean value $\mean{w}$, we obtain
\begin{align}\label{delta:Apri:1}
\norm{W(u)}_{L^{\infty}(0,T;L^{1})} + \norm{u}_{L^{\infty}(0,T;H^{1})} + \norm{w}_{L^{2}(0,T;H^{1})} \leq C.
\end{align}
Next, taking the time derivative of \eqref{CHd:2}, testing with $w$, and adding to the resulting equation the equality obtained from testing \eqref{CHd:1} with $\pd_{t} u$ gives
\begin{align*}
\norm{\pd_{t}u}_{L^{2}}^{2} + \frac{\dd}{\dt} \frac{1}{2} \norm{w}_{L^{2}}^{2} = \int_{\Omega} W''(u) \, \pd_{t} u \, w + \lambda (I - u) \, \pd_{t} u \dx.
\end{align*}
Applying Young's inequality, the fact that $W''$ is bounded, and Gronwall's inequality leads to
\begin{align}\label{delta:Apri:3}
\norm{w}_{L^{\infty}(0,T;L^{2})} + \norm{\pd_{t} u}_{L^{2}(0,T;L^{2})} \leq C.
\end{align}
Let us point out that it suffices to consider initial conditions $u_{0} \in H^{2}(\Omega)$, so that $w(0) := -\Laplace u_{0} + W'(u_{0}) \in L^{2}(\Omega)$.  Then, as $\pd_{t}u - \lambda(I - u) \in L^{2}(0,T;L^{2}(\Omega))$ and $w - W'(u) \in L^{\infty}(0,T;L^{2}(\Omega))$, by elliptic regularity theory we deduce that
\begin{align}\label{delta:Apri:4}
\norm{w}_{L^{2}(0,T;H^{2})} + \norm{u}_{L^{\infty}(0,T;H^{2})} \leq C.
\end{align}
The estimates \eqref{delta:Apri:1}-\eqref{delta:Apri:4} ensure that the solutions to \eqref{CHd:1}-\eqref{CHd:2} are sufficiently regular in order to obtain the a priori estimates \eqref{Apri:1}, \eqref{Apri:2}, \eqref{Apri:3}, \eqref{Apri:4}, \eqref{Apri:5}, and \eqref{Apri:6}.

\section{The stationary problem}\label{sec:stat}

The stationary problem of \eqref{Vineq} is given as
\begin{subequations}\label{stat:prob}
\begin{alignat}{3}
- \Laplace w = \lambda (I - u) & \text{ in } \Omega, \label{stat:1} \\
- \left ( w + \eps^{-1} u, v - u \right ) + (\eps \nabla u, \nabla (v - u)) \geq 0 & \text{ for all } v \in \KK, \label{stat:2}  \\
\pdnu u = \pdnu w = 0 & \text{ on } \pd \Omega.
\end{alignat}
\end{subequations}
We point out that due to the boundary condition for $w$, upon integrating, \eqref{stat:1} leads to the condition
\begin{align}\label{stat:mean}
\int_{\Omega} \lambda (I - u) \dx = \alpha \int_{\Omega \setminus D}  (I - u) \dx = 0.
\end{align}
In particular, the spatial mean of $u$ in the undamaged region is equal to the spatial mean of the image data $I$ in the undamaged region.  Let us introduce the subspaces
\begin{align*}
V' & := \{ g \in (H^{1}(\Omega))' : \inner{g}{1}_{H^{1}} = 0 \} \subset (H^{1}(\Omega))', \\
L^{2}_{0}(\Omega) & := \{ h \in L^{2}(\Omega) : \mean{h} = 0 \} \subset L^{2}(\Omega)
\end{align*}
and the operator $\NN : V' \to H^{1}(\Omega) \cap L^{2}_{0}(\Omega)$ as the solution operator to the Neumann--Laplacian:
\begin{align}\label{defn:N}
\NN(f) = v \quad \Longleftrightarrow \quad \int_{\Omega} \nabla v \cdot \nabla \zeta \dx = \inner{f}{\zeta}_{H^{1}} \quad \text{ for all } \zeta \in H^{1}(\Omega).
\end{align}
Then, due to \eqref{stat:mean}, we can express $w$ as a sum of its spatial mean $\mean{w}$ and its mean-free part $\NN \left ( \lambda (I - u ) \right )$, i.e., $w= \mean{w} + \NN(\lambda(I - u))$, and this allows us to express \eqref{stat:prob} as one variational inequality
\begin{align}\label{stat:one:vineq}
- \left (\mean{w} + \eps^{-1} u + \NN(\lambda (I - u)), v - u \right ) + \left ( \eps \nabla u, \nabla (v - u) \right ) \geq 0 \text{ for all } v \in \KK.
\end{align}  
In the classical Cahn--Hilliard situation, i.e., $\lambda=0$ in
\eqref{stat:prob}, the existence of a solution can be easily shown using variational arguments.  As \eqref{stat:prob} is not known to be the first variation of a functional we need to derive another approach.  Our result concerning the existence of a solution $(u,w)$ to the stationary problem \eqref{stat:prob} is formulated as follows.

\begin{thm}\label{thm:2}
Under assumptions \eqref{ass:Omega} and \eqref{ass:I}, for every $\alpha, \eps > 0$, there exists at least one pair $(u,w) \in H^{2}_{n}(\Omega) \times (H^{2}_{n}(\Omega) \cap W^{2,p}(\Omega))$ for any $p < \infty$ satisfying $\abs{u} \leq 1$ a.e. in $\Omega$, $\int_{\Omega \setminus D} u - I \dx = 0$, and
\begin{subequations}\label{stat:weakform}
\begin{alignat}{3}
\Laplace w + \lambda (I - u)  = 0 & \text{ a.e. in } \Omega,  \label{stat:weak:1} \\
-(w + \eps^{-1} u, v - u) + (\eps \nabla u, \nabla (v - u))  \geq 0 & \text{ for all } v \in \KK. \label{stat:weak:2}
\end{alignat}
\end{subequations}
\end{thm} 

The closest companion result to Theorem~\ref{thm:2} is the existence result sketched in \cite[Appendix~A]{BHS} for Neumann boundary conditions.  We point out that the system studied in \cite{BHS} is
\begin{align*}
-\eps \Laplace u + \eps^{-1} W_{\mathrm{qu}}'(u) = \NN \left ( \lambda (I - u) - \mean{\lambda (I - u)} \right ) & \text{ in } \Omega, \\
\pdnu u = 0 & \text{ on } \pd \Omega,
\end{align*}
and compared to \eqref{stat:one:vineq} it appears that the spatial mean $\mean{w}$ has been neglected and the fidelity term $\lambda (I - u)$ has been modified to have zero spatial mean.

\subsection{Approximation scheme}

Let $g \in C^{\infty}(\R)$ such that $0 \leq g(s) \leq 1$ for all $s \in \R$ and
\begin{align}\label{g}
g(s) = \begin{cases}
1 & \text{ if } s \geq 3, \\
0 & \text{ if } s \leq 2.
\end{cases}
\end{align}
Then, we define for $v \in L^{2}(\Omega)$
\begin{align}
\label{defn:F}
F(v) := C_{F} g  \left ( \norm{v}_{L^{2}}^{2} / \abs{\Omega} \right ),
\end{align}
where $C_{F}$ is a positive constant to be specified later.  We consider a sequence of solutions $\{u_{\delta}\}_{\delta > 0} \subset W := H^{2}_{n}(\Omega)$, indexed by $\delta > 0$, such that 
\begin{subequations}
\begin{alignat}{3}
F(u_{\delta}) u_{\delta} + \sqrt{\delta} \beta_{\delta}(u_{\delta}) + \eps \Laplace^{2} u_{\delta} - \eps^{-1} \Laplace \beta_{\delta}(u_{\delta}) & = \lambda (I - u_{\delta}) - \eps^{-1} \Laplace u_{\delta} && \text{ in } \Omega, \\
\pdnu u_{\delta} = \pdnu \Laplace u_{\delta} & = 0 && \text{ on } \pd \Omega,
\end{alignat}
\end{subequations}
which holds in the following weak sense:
\begin{equation}\label{stat:approx}
\begin{aligned}
& \int_{\Omega} F(u_{\delta})u_{\delta} \zeta + \sqrt{\delta} \beta_{\delta}(u_{\delta}) \zeta + \eps \Laplace u_{\delta} \Laplace \zeta + \eps^{-1} \nabla \beta_{\delta}(u_{\delta}) \cdot \nabla \zeta \dx \\
& \quad = \int_{\Omega} \lambda (I - u_{\delta}) \zeta + \eps^{-1} \nabla u_{\delta} \cdot \nabla \zeta \dx
\end{aligned}
\end{equation}
for all $\zeta \in W$.  For $i = 1,2$, we define operators $A_{i} :W \to W'$ as
\begin{equation}
\begin{aligned}
\inner{A_{1}u}{\zeta}_{W} & := \int_{\Omega} \sqrt{\delta} \beta_{\delta}(u) \zeta + \eps \Laplace u \Laplace \zeta, \\
\inner{A_{2}u}{\zeta}_{W} & := \int_{\Omega} \eps^{-1} \nabla \beta_{\delta}(u) \cdot \nabla \zeta - \lambda (I - u) \zeta - \eps^{-1}  \nabla u \cdot \nabla \zeta + F(u)u \zeta \dx,
\end{aligned}
\end{equation}
so that \eqref{stat:approx} is equivalent to $\inner{A_{1}u_{\delta} + A_{2} u_{\delta}}{\zeta}_{W} = 0$ for all $\zeta \in W$.  Note that the operator $A_{1}$ is monotone and hemicontinuous (see \cite[\S~26.1]{Zeidler}), whereas the operator $A_{2}$ is strongly continuous, i.e., $u_{n} \rightharpoonup u$ in $W$ implies $A_{2} u_{n} \to A_{2} u$ in $W'$.  This follows from the continuity and sublinear growth of $\beta_{\delta}$, the continuity and boundedness of $F$, the compact embedding $W \subset \subset H^{1}(\Omega)$, and Rellich's theorem.  Now, the application of  \cite[Theorem~27.6]{Zeidler} yields that $A= A_{1} + A_{2}$ is a pseudomonotone operator.  We further claim that $A = A_{1} + A_{2}$ is coercive, i.e.,
\begin{align*}
\underset{\norm{u}_{W} \to \infty}{\lim} \frac{\inner{Au}{u}_{W}}{\norm{u}_{W}} = \infty.
\end{align*}
This follows from the identity
\begin{align*}
\inner{Au}{u}_{W} & = \int_{\Omega} F(u) \abs{u}^{2} + \sqrt{\delta} \beta_{\delta}(u) u + \eps \abs{\Laplace u}^{2} + \eps^{-1} \beta_{\delta}'(u) \abs{\nabla u}^{2} \dx \\
& \quad  - \int_{\Omega}  \lambda (I - u)u + \eps^{-1} \abs{\nabla u}^{2} \dx.
\end{align*}
By the property \eqref{g} of the smooth function $g$, for $\norm{u}_{L^{2}}^{2} \geq 3 \abs{\Omega}$ we have
\begin{align*}
\int_{\Omega} F(u) \abs{u}^{2} \dx \geq C_{F} \norm{u}_{L^{2}}^{2}.
\end{align*}
Furthermore we have
\begin{align*}
-\int_{\Omega} \lambda (I-u)u \dx = \int_{\Omega} \lambda (\abs{u}^{2} - Iu) \dx \geq \int_{\Omega} \frac{\lambda}{2} \abs{u}^{2} \dx - \frac{\alpha}{2} \abs{\Omega \setminus D}.
\end{align*}
Using inequality \eqref{inequDelta} with $\eta = \frac{\eps^{2}}{2}$, the monotonicity of $\beta_{\delta}$ together with the last two inequalities we obtain for $\norm{u}_{L^{2}}^{2} \geq 3 \abs{\Omega}$,
\begin{align*}
\inner{Au}{u}_{W} & \geq \int_{\Omega} (C_{F} - \tfrac{1}{2} \eps^{-3}) \abs{u}^{2} + \frac{1}{2}\eps \abs{\Laplace u}^{2} \dx  - \frac{\alpha}{2} \abs{\Omega \setminus D}.
\end{align*}
Choosing $C_{F} = \eps^{-3}$ and using the fact that $\norm{f}^{2} := \int_{\Omega} \abs{\Laplace f}^{2} + \abs{f}^{2} \dx$ is equivalent to the $W$-norm we obtain coercivity of $A$.

For each $\delta \in (0,1)$, by \cite[Theorem~27.A]{Zeidler} there exists a solution $u_{\delta} \in W$ to the abstract equation $A u_{\delta} = 0$.  We now define
\begin{align*}
w_{\delta} = - \eps \Laplace u_{\delta} + \eps^{-1}(\beta_{\delta}(u_{\delta}) - u_{\delta}).
\end{align*}
The equality $A u_{\delta} =0$ implies that for all $\zeta \in W$,
\begin{align*}
\int_{\Omega} w_{\delta} \Laplace \zeta \dx = \int_{\Omega} f \zeta \dx
\end{align*}
holds, where
\begin{align*}
f = F(u_{\delta})u_{\delta} + \sqrt{\delta} \beta_{\delta}(u_{\delta}) - \lambda(I-u_{\delta})\in L^{2}_{0}(\Omega).
\end{align*}
Indeed, testing with $\zeta=1$ in \eqref{stat:approx}	leads to $\int_{\Omega}  f \dx =0$.

We claim that $w_{\delta} \in W$.  Indeed, there exists a weak solution $\hat{w} \in H^{1}(\Omega)$ to the variational problem
\begin{align}\label{eq2}
\int_{\Omega} \nabla \hat{w} \cdot \nabla \eta \dx = \int_{\Omega} -f \eta \dx
\end{align}
for all $\eta \in H^{1}(\Omega)$.  Since the solution to \eqref{eq2} is uniquely determined up to a constant, we choose $\hat{w}$ such that $\int_{\Omega} \hat{w} \dx = \int_{\Omega} w_{\delta} \dx$.  Then, elliptic regularity theory gives that the weak solution $\hat{w}$ also fulfills $\hat{w} \in W$.  The
function $z= w_{\delta} - \hat{w}$ in turn satisfies
\begin{align}\label{eq3}
\int_{\Omega} z \Laplace \zeta \dx = 0  \text{ for all } \zeta \in W.
\end{align}
We now solve the auxiliary problem
\begin{align*}
\Laplace y = z \text{ in } \Omega, \quad \pdnu y = 0 \text{ on } \pd \Omega \quad \text{ with } \mean{y} = 0,
\end{align*}
which is possible as $\int_{\Omega} z \dx = 0$.  Standard elliptic theory yields that $y \in W$, and choosing $\zeta = y$ in
\eqref{eq3} now gives $\norm{z}_{L^{2}}^{2} = 0$.  This implies that $z=0$ and $w_{\delta} = \hat{w} \in W$.

\subsection{Uniform estimates}

The pair $(u_{\delta}, w_{\delta}) \in W \times W$ fulfills the following weak
formulation: For all $\zeta \in W$ it holds that
\begin{subequations}
\begin{alignat}{3}
\int_{\Omega} w_{\delta} \zeta \dx = \int_{\Omega} \eps \nabla u_{\delta} \cdot \nabla \zeta + \eps^{-1} (\beta_{\delta}(u_{\delta}) - u_{\delta})\zeta \dx, \label{stat:delta:1} \\
\int_{\Omega} F(u_{\delta}) u_{\delta} \zeta + \sqrt{\delta} \beta_{\delta}(u_{\delta}) \zeta + \nabla w_{\delta} \cdot \nabla \zeta \dx = \int_{\Omega} \lambda (I - u_{\delta}) \zeta \dx. \label{stat:delta:2}
\end{alignat}
\end{subequations}
As $W$ is dense in $H^{1}(\Omega)$, these equalities also hold for all $\zeta \in H^{1}(\Omega)$.

\paragraph{First estimate.} Similar to the previous argument in the proof of the coercivity of $A$, we have
\begin{align*}
& \int_{\Omega}  F(u_{\delta}) \abs{u_{\delta}}^{2} + \sqrt{\delta} \beta_{\delta}(u_{\delta}) u_{\delta} + \eps \abs{\Laplace u_{\delta}}^{2} + \eps^{-1} \beta_{\delta}'(u_{\delta}) \abs{\nabla u_{\delta}}^{2}  - \eps^{-1} \abs{\nabla u_{\delta}}^{2} \dx  \\
& \quad = \int_{\Omega} \lambda (I - u_{\delta})u_{\delta} \dx \leq - \frac{1}{2} \int_{\Omega} \lambda \abs{u_{\delta}}^{2} + \frac{\alpha}{2} \abs{\Omega \setminus D}.
\end{align*}
Using \eqref{inequDelta} with $\eta = \frac{\eps}{2}$, we have
\begin{align*}
- \eps^{-1} \norm{\nabla u_{\delta}}_{L^{2}}^{2} \geq - \frac{\eps}{2} \norm{\Laplace u_{\delta}}_{L^{2}}^{2} - \frac{1}{2} \eps^{-3} \norm{u_{\delta}}_{L^{2}}^{2},
\end{align*}
and this together with the nonnegativity of $\lambda$ and the monotonicity of $\beta_{\delta}$ gives
\begin{align*}
\int_{\Omega} \sqrt{\delta} \beta_{\delta}(u_{\delta})u_{\delta} + (F(u_{\delta})-\tfrac{1}{2} \eps^{-3}) \abs{u_{\delta}}^{2} + \frac{\eps}{2} \norm{\Laplace u_{\delta}}_{L^{2}}^{2} \leq \frac{\alpha}{2} \abs{\Omega \setminus D}.
\end{align*}
If $\norm{u_{\delta}}_{L^{2}}^{2} \geq 3\abs{\Omega}$ we have for $C_F = \eps^{-3}$
\begin{align*}
\int_{\Omega} \sqrt{\delta} \beta_{\delta}(u_{\delta})u_{\delta} \dx + \eps^{-3} \norm{u_{\delta}}_{L^{2}}^{2} + \eps \norm{\Laplace u_{\delta}}_{L^{2}}^{2} \leq \alpha \abs{\Omega \setminus D}.
\end{align*}
If $\norm{u_{\delta}}_{L^{2}}^{2} \leq 3 \abs{\Omega}$, neglecting the nonnegative term $F(u_{\delta}) \abs{u_{\delta}}^2$, we have
\begin{align*}
& \int_{\Omega} \sqrt{\delta} \beta_{\delta}(u_{\delta})u_{\delta} \dx + \norm{u_{\delta}}_{L^{2}}^{2} + \frac{\eps}{2} \norm{\Laplace u_{\delta}}_{L^{2}}^{2} \\
& \quad \leq 3 \abs{\Omega} + \frac{\alpha}{2} \abs{\Omega \setminus D} + \frac{1}{2} \eps^{-3} \norm{u_{\delta}}_{L^{2}}^{2} \leq 3 \abs{\Omega} ( 1 + \tfrac{1}{2} \eps^{-3}) + \frac{\alpha}{2} \abs{\Omega \setminus D}.
\end{align*}
This implies that $\{u_{\delta}\}_{\delta \in (0,1)}$ is bounded in $W$.  

\paragraph{Second estimate.} Next, we observe
\begin{align*}
(\beta_{\delta}(z)z)' = \beta_{\delta}'(z)z + \beta_{\delta}(z) = \beta_{\delta}'(z)z + \hat{\beta}_{\delta}'(z),
\end{align*}
where we recall the antiderivative $\hat{\beta}_{\delta}$ of $\beta_{\delta}$.  As $\beta_{\delta}'(z) \geq 0$, we see that
\begin{align*}
(\beta_{\delta}(z)z)' \geq \hat{\beta}_{\delta}'(z) & \text{ for } z > 0, \\
(\beta_{\delta}(z)z)' \leq \hat{\beta}_{\delta}'(z) & \text{ for } z < 0.
\end{align*}
Together with the fact that $\beta_{\delta}(0) = 0 = \hat{\beta}_{\delta}(0)$ we obtain
\begin{align*}
\hat{\beta}_{\delta}(z) \leq \beta_{\delta}(z)z  \text{ for all } z \in \R.
\end{align*}
Thanks to the uniform bound $\int_{\Omega} \sqrt{\delta} \beta_{\delta}(u_{\delta})u_{\delta} \dx \leq C$, we obtain
\begin{align}\label{stat:hat:beta}
\sqrt{\delta} \int_{\Omega} \hat{\beta}_{\delta}(u_{\delta}) \dx \leq C.
\end{align}
From the definition \eqref{hat:beta:approx} we have $\hat{\beta}_{\delta}(u_{\delta}) = \delta^{-1} \psi(u_{\delta})$, where 
\begin{align*}
\psi(s):= \begin{cases}
\frac{1}{2} (s-1)^{2} & \text{ if } s \geq 1, \\
0 & \text{ if } s\in [-1,1],\\
\frac{1}{2} (s+1)^{2} & \text{ if } s \leq -1.
\end{cases}
\end{align*}
Hence, from \eqref{stat:hat:beta} we infer
\begin{align*}
\int_{\Omega} \psi(u_{\delta}) \leq C\sqrt{\delta} \to 0 \text{ as } \delta \to 0.
\end{align*}
Thanks to the uniform boundedness of $u_{\delta}$ in $H^{2}(\Omega)$, for a  nonrelabelled subsequence
\begin{align*}
u_{\delta} & \rightharpoonup u \text{ in } H^{2}(\Omega), \\
u_{\delta} & \to u \text{ in } H^{1}(\Omega) \text{ and a.e. in } \Omega.
\end{align*}
Continuity of $\psi$ then yields $\psi(u_{\delta}) \to \psi(u)$ a.e. in $\Omega$, and as a result of Fatou's lemma we infer 
\begin{align*}
\psi(u) = 0 \text{ a.e. in } \Omega \quad \Longrightarrow \quad u \in [-1,1] \text{ a.e. in } \Omega.
\end{align*}
In turn, this bound gives $\norm{u}_{L^{2}}^{2} \leq \abs{\Omega}$, and as $\norm{u_{\delta}}_{L^{2}}^{2} \to \norm{u}_{L^{2}}^{2}$, by the property \eqref{g} of the smooth function $g$, we see that the term $F(u_{\delta})u_{\delta}$ in \eqref{stat:delta:2} is zero for $\delta$ sufficiently small. 

\paragraph{Third estimate.} We test \eqref{stat:delta:1} with $-\Delta w_\delta$ and obtain
\begin{align*}
\norm{\nabla w_{\delta}}_{L^{2}}^{2} = \int_{\Omega} \eps \Laplace u_{\delta} \Laplace w_{\delta} - \eps^{-1}(\beta_{\delta}(u_{\delta}) - u_{\delta})) \Laplace w_{\delta} \dx.
\end{align*}
We then test \eqref{stat:delta:2} with $ \eps^{-1} \beta_\delta(u_\delta)$ and $- \eps \Delta
u_\delta$ and obtain for sufficiently small $\delta$
\begin{align*}
\sqrt{\delta} \eps^{-1} \norm{\beta_{\delta}(u_{\delta})}_{L^{2}}^{2} + \int_{\Omega} \eps^{-1} \lambda (u_{\delta} - I) \beta_{\delta}(u_{\delta}) \dx = \int_{\Omega} \eps^{-1} \beta_{\delta}(u_{\delta}) \Laplace w_{\delta} \dx, \\
\int_{\Omega} \eps \sqrt{\delta} \beta_{\delta}'(u_{\delta}) \abs{\nabla u_{\delta}}^{2} \dx = -\int_{\Omega} \eps \Laplace w_{\delta} \Laplace u_{\delta} + \eps \lambda (I - u_{\delta}) \Laplace u_{\delta} \dx.
\end{align*}
Adding the last three equalities gives
\begin{equation}\label{stat:sqr:beta}
\begin{aligned}
& \norm{\nabla w_{\delta}}_{L^{2}}^{2} + \sqrt{\delta} \eps^{-1} \norm{\beta_{\delta}(u_{\delta})}_{L^{2}}^{2} +  \int_{\Omega} \eps \sqrt{\delta} \beta_{\delta}'(u_{\delta}) \abs{\nabla u_{\delta}}^{2} + \eps^{-1} \lambda \underbrace{(u_{\delta} - I) \beta_{\delta}(u_{\delta})}_{\geq 0 \text{ by } \eqref{I-u:beta}} \dx \\
& \quad = \int_{\Omega} \eps^{-1} u_{\delta} \Laplace w_{\delta} - \eps \lambda (I - u_{\delta}) \Laplace u_{\delta} \dx  \leq \frac{1}{2} \norm{\nabla w_{\delta}}_{L^{2}}^{2} + C \left ( 1 + \norm{u_{\delta}}_{H^{2}}^{2} \right ),
\end{aligned}
\end{equation}
where $C$ is a positive constant depending on $\eps$ and $\alpha$ but not on $\delta$.  Hence, we obtain that $\{\nabla w_{\delta}\}_{\delta \in (0,1)}$ is bounded in $L^{2}(\Omega)$.  Furthermore, using $\sqrt{\delta} \norm{\beta_{\delta}(u_{\delta})}_{L^{2}}^{2} \leq C$, we see that
\begin{align*}
\norm{\sqrt{\delta} \beta_{\delta}(u_{\delta})}_{L^{2}}^{2} = \delta \norm{\beta_{\delta}(u_{\delta})}_{L^{2}}^{2} \leq C \sqrt{\delta} \to 0 \text{ as } \delta \to 0.
\end{align*}
This implies that $\sqrt{\delta} \beta_{\delta}(u_{\delta}) \to 0$ as $\delta \to 0$, and together with $F(u_{\delta}) = 0$ for $\delta$ sufficiently small, we obtain from \eqref{stat:delta:2} with $\zeta = 1$ and the strong convergence of $u_{\delta}$ to $u$ the identity
\begin{align*}
\int_{\Omega} \lambda (I - u) \dx = 0 \quad \Longrightarrow \quad \int_{\Omega \setminus D} u \dx = \int_{\Omega \setminus D} I \dx.
\end{align*}
From the above we can also infer that the mean value of $u$ lies in the interval $(-1,1)$.  Indeed, thanks to \eqref{ass:I}, the image function $I$ is not identically equal to $1$ or $-1$ in $\Omega \setminus D$, and so
\begin{align}\label{stat:mean:limit}
\abs{\frac{1}{\abs{\Omega}} \int_{\Omega} u \dx} = \frac{1}{\abs{\Omega}} \abs{\int_{\Omega \setminus D} I \dx + \int_{D} u \dx} < \frac{1}{\abs{\Omega}} (\abs{\Omega \setminus D} + \abs{D}) = 1.
\end{align}

\paragraph{Fourth estimate.} To obtain uniform estimates on the mean values of $w_{\delta}$, we argue as in the time-dependent case.  Testing \eqref{stat:delta:1} with $\zeta = \pm 1$ gives
\begin{align*}
\abs{\int_{\Omega} w_{\delta} \dx} \leq \eps^{-1} \int_{\Omega} \abs{\beta_{\delta}(u_{\delta})} + \abs{u_{\delta}} \dx.
\end{align*}
Meanwhile, testing \eqref{stat:delta:1} with $u_{\delta}$ leads to
\begin{align*}
\eps \norm{\nabla u_{\delta}}_{L^{2}}^{2} + \int_{\Omega} \eps^{-1} u_{\delta} \beta_{\delta}(u_{\delta}) \dx = \int_{\Omega} w_{\delta} u_{\delta} + \eps^{-1} \abs{u_{\delta}}^{2} \dx.
\end{align*}
Then, using the fact that $\abs{\beta_{\delta}(r)} \leq r \beta_{\delta}(r)$ for all $r \in \R$, we combine the above two estimates to obtain
\begin{align}\label{stat:mean:mu:1}
\abs{\int_{\Omega} w_{\delta} \dx} \leq C \eps^{-1} \left ( 1+ \norm{u_{\delta}}_{L^{2}}^{2}\right ) + \int_{\Omega} w_{\delta} u_{\delta} \dx.
\end{align}
Next, recalling the operator $\mathcal{N} : V' \to H^{1}(\Omega) \cap L^{2}_{0}(\Omega)$ defined in \eqref{defn:N}, we now test \eqref{stat:delta:2} with $\mathcal{N}(u_{\delta} - \mean{u_{\delta}})$ (with $\delta$ sufficiently small so that $F(u_{\delta}) = 0$) to arrive at
\begin{equation}\label{stat:mean:mu:2}
\begin{aligned}
& \int_{\Omega} (\lambda (I - u_{\delta}) - \sqrt{\delta} \beta_{\delta}(u_{\delta})) \mathcal{N}(u_{\delta} - \mean{u_{\delta}}) \dx = \int_{\Omega} \nabla w_{\delta} \cdot \nabla \mathcal{N}(u_{\delta} - \mean{u_{\delta}}) \dx \\
& \quad  = \int_{\Omega} w_{\delta} (u_{\delta} - \mean{u_{\delta}}) \dx = \int_{\Omega} w_{\delta} u_{\delta} \dx - \mean{u_{\delta}} \int_{\Omega} w_{\delta} \dx.
\end{aligned}
\end{equation}
Employing the estimate
\begin{align*}
\norm{\mathcal{N}(u_{\delta} - \mean{u_{\delta}})}_{H^{1}} \leq c \norm{u_{\delta} - \mean{u_{\delta}}}_{L^{2}} \leq c^{2} \norm{\nabla u_{\delta}}_{L^{2}},
\end{align*}
where $c$ is the positive constant from the Poincar\'{e} inequality, we obtain from \eqref{stat:mean:mu:1}-\eqref{stat:mean:mu:2}
\begin{equation}\label{stat:mean:mu:3}
\begin{aligned}
(1 - \abs{\mean{u_{\delta}}}) \abs{\int_{\Omega} w_{\delta} \dx} & \leq C\eps^{-1} \left ( 1 + \norm{u_{\delta}}_{L^{2}}^{2} \right ) \\
& \quad + C \left ( 1 + \norm{u_{\delta}}_{L^{2}} + \sqrt{\delta} \norm{\beta_{\delta}(u_{\delta})}_{L^{2}} \right ) \norm{\nabla u_{\delta}}_{L^{2}}.
\end{aligned}
\end{equation}
By virtue of the strong convergence of $u_{\delta}$ to $u$ in $L^{1}(\Omega)$, as well as the property \eqref{stat:mean:limit} that $\mean{u} \in (-1,1)$, there exists a $\delta_{0} > 0$ such that for all $\delta < \delta_{0}$, it holds that $\abs{\mean{u_{\delta}}} < 1$, and so the prefactor $(1-\abs{\mean{u_{\delta}}})$ on the left-hand side of \eqref{stat:mean:mu:3} is uniformly bounded away from zero for all $\delta < \delta_{0}$.  Furthermore, by \eqref{stat:sqr:beta}, we have $\sqrt{\delta} \norm{\beta_{\delta}(u_{\delta})}_{L^{2}}^{2} \leq C$, and so the right-hand side of \eqref{stat:mean:mu:3} is uniformly bounded.  Hence, we infer that $\{\mean{w_{\delta}}\}_{\delta \in (0,\delta_{0})}$ is uniformly bounded.  Together with \eqref{stat:mean:mu:3}, the application of the Poincar\'{e} inequality yields the uniform boundedness of $w_{\delta}$ in $H^{1}(\Omega)$.  Thus, along a nonrelabelled subsequence
\begin{align*}
w_{\delta} & \rightharpoonup w \text{ in } H^{1}(\Omega), \\
w_{\delta} & \to w \text{ in } L^{2}(\Omega) \text{ and a.e. in } \Omega.
\end{align*}
By testing \eqref{stat:delta:1} with $\beta_{\delta}(u_{\delta})$, we obtain as in the time-dependent case
\begin{align*}
\norm{\beta_{\delta}(u_{\delta})}_{L^{2}}^{2} \leq C.
\end{align*}
Using the fact that $\sqrt{\delta} \beta_{\delta}(u_{\delta}) \to 0$ in $L^{2}(\Omega)$ we can pass to the limit in \eqref{stat:delta:2} to obtain the equality
\begin{align*}
(\nabla w, \nabla \zeta) - (\lambda (I - u), \zeta) & = 0 \quad \text{ for all } \zeta \in H^{1}(\Omega).
\end{align*}
Using the fact that $\lambda (I - u) \in L^{\infty}(\Omega)$ and elliptic regularity theory gives $w \in W^{2,p}(\Omega)$ for all $p < \infty$.  This allows us to express the above variational equality as \eqref{stat:weak:1}.  Meanwhile, passing to the limit in \eqref{stat:delta:1} and arguing as in the time-dependent case leads to \eqref{stat:weak:2}.  This concludes the proof of Theorem~\ref{thm:2}.

\section{Numerical implementation}\label{sec:num}
In this section we first derive a finite element approximation of \eqref{weakform}, and then we display some numerical results obtained using this approximation.  Let us mention that there are two methods to solve the time-dependent problem \eqref{weakform}: The first is to propose a discretization of the Moreau--Yosida approximation \eqref{CHdelta} for fixed $\delta > 0$, and the second is to solve the variational inequality \eqref{weakform} directly.  Note that for fixed $\delta > 0$, \eqref{CHdelta} is essentially the usual modified Cahn--Hilliard equation with a $C^{2}$-potential, and one can apply finite elements or spectral methods for fast inpainting.  One issue of this approach is the choice of $\delta \ll 1$ relative to other parameters $\eps$ and $\alpha$.  It is desirable to have $\delta$ rather small, and it is recommended in \cite{Bosch} to focus on a discretization based on finite elements as opposed to spectral methods, due to a significantly higher number of iterations needed for a fast Fourier-transform (FFT) based implementation compared to a standard finite element implementation, and an undesirable increase in FFT iterations as the penalization parameter $\delta$ decreases.  Furthermore, the finite element framework is favoured by \cite{Bosch} over the finite difference framework as the former allows the authors to compute for missing information on arbitrary domains, such as complex three-dimensional objects.

In contrast, solving the variational inequality \eqref{weakform} directly, as we do below, avoids the issues involving the parameter $\delta$.  Ideally, small values of the parameter $\eps$ should be used during implementation for a closer approximation to the original binary image, but this requires more resolution in the interfacial regions $\{ -1 < u < 1 \}$ separating bulk regions $\{ u = \pm 1 \}$.  For implementations with finite differences, adaptivity and especially local refinement cannot be done as easily and as flexibly as, for example, in the context of finite element methods.
Therefore, we choose to employ a finite element discretization of the variational inequality \eqref{weakform}, due to the well-established literature in this area and efficient solvers for Cahn--Hilliard variational inequalities, which in particular also allow for adaptivity.  We mention that there are also error analyses for space-semidiscrete and fully discrete finite element schemes for \eqref{CH} with Dirichlet boundary conditions; see \cite{Poole} for more details.

\subsection{Finite element scheme}
Let $\mathcal{T}$ be a regular triangulation of $\Omega$ into disjoint open simplices; associated with $\mathcal{T}$ is the piecewise linear finite element space
\begin{align*}
S_{h} :=  \left \{ \varphi \in C^{0}(\overline{\Omega}) 
\, : \,  \varphi \vert_{T} \in P_{1}(T) \; \text{ for all }~T \in \mathcal{T} \right \} \subset H^{1}(\Omega),
\end{align*}
where we denote by $P_{1}(T)$ the set of all affine linear functions on $T$.  We now introduce a finite element approximation of \eqref{weakform}; see \cite{BNS} for more details.

Given $u_{h}^{n-1} \in K_{h} := \{ \chi \in S_{h} \, : \, \abs{\chi} \leq 1 \}$ find $\{u_{h}^{n},\, w_{h}^{n}\} \in K_{h} \times S_{h}$ such that for all $\eta_{h} \in S_{h}$ and $\zeta_{h} \in K_{h}$,
\begin{subequations}\label{FE}
\begin{align}
\frac{1}{\tau} ( u_{h}^{n} - u_{h}^{n-1}, \eta_{h})_{h} + ( \nabla w_{h}^{n}, \nabla \eta_{h}) & = (\lambda(x)(I -u_{h}^{n-1}), \eta_{h})_{h}, \label{FE:u} \\
\left ( w_{h}^{n} + \frac{1}{\eps} u_{h}^{n-1} , \zeta_{h} - u_{h}^{n} \right )_{h} & \leq  \eps (\nabla u_{h}^{n}, \nabla (\zeta_{h} - u_{h}^{n})), \label{FE:w} 
\end{align}
\end{subequations}
where $\tau$ denotes the time step, and $(\eta_{1},\eta_{2})_{h} := \int_{\Omega}\pi_{h}(\eta_{1}(x) \eta_{2}(x)) \dx$, where on each triangle $\pi_{h}(\eta_1 \eta_2)$ is taken to be an affine interpolation of the values of $\eta_{1} \eta_{2}$ at the nodes of the triangle.

As a stopping criteria for the scheme we follow the procedure of \cite{BEG,Bosch}, in which the inpainted image is computed in a two step process.  In the first step the scheme is solved with a relatively large value of $\eps$ denoted by $\eps_{1}$, until the stopping criteria of 
\begin{align}\label{eq:sc}
\norm{u_{h}^{n} - u_{h}^{n-1}}^2 \leq 5.0\cdot 10^{-6}
\end{align}
is fulfilled.  To sharpen the edges of the approximate solution we reduce the value of $\eps$ to a smaller value $\eps_{2}$ and increase the value of $\alpha$ to a large value $\alpha_{2}$, and the computation is terminated when the stopping criteria \eqref{eq:sc} is reached with these new parameters $\eps_2$ and $\alpha_2$.  We denote the solution of the terminated computation by $\tilde{u}_{h}$.

We use the finite element toolbox ALBERTA 2.0 \cite{alberta} for adaptivity, and we implemented a similar mesh refinement strategy to that in \cite{BNS}, i.e., a fine mesh is constructed where $\abs{u_{h}^{n-1}} < 1.0$ with a coarser mesh present in the bulk regions $\abs{u_{h}^{n-1}} =1.0$.  In the simulations in the next section, unless otherwise stated we set the minimal diameter of an element $h_{min} = 3.9 \cdot 10^{-3}$ and the maximal diameter $h_{max} = 2.2 \cdot 10^{-2}$.

\subsection{Numerical simulations}

We conclude with some numerical simulations; first we apply the finite element approximation \eqref{FE:u}-\eqref{FE:w} to binary images and then we consider grayscale images by using the generalization of \eqref{CH:1}-\eqref{CH:2} presented in \cite{sb}. In order to cancel out the smoothing effects and to focus on the reconstruction of the damaged region, rather than plotting the solution $\tilde{u}_h$ we instead plot $P(\tilde{u}_{h}):=\chi_{\{\tilde{u}_h\geq 0\}}-\chi_{\{\tilde{u}_h<0\}}$ which is the projection of $\tilde{u}_{h}$ into the set of binary images. The values of the parameters that we used in the simulations for Figures \ref{f:comp} - \ref{f:gs3} are displayed in Table \ref{t:param}.

\begin{table}[!h]\begin{center}
 \begin{tabular}{ |c||c|c|c|c|c|c| }
 \hline
~ & $\eps_1$ & $\eps_2$ & $\alpha$ & $\alpha_2$ & $\tau$ \\
 \hline
 \hline
 Figure \ref{f:comp} &$0.04$ & $0.00\dot{3}$ & $8.0\cdot 10^{3}$ & $1.0\cdot10^{5}$ & $1.0\cdot10^{-5}$\\
Figure \ref{f:1} &$0.04$ & $0.00\dot{3}$ & $8.0\cdot 10^{3}$ & $1.0\cdot10^{5}$ & $1.0\cdot10^{-5}$\\
Figure \ref{f:2} &$0.0125$ & $0.00\dot{3}$ & $1.0\cdot10^{6}$& $3.0\cdot10^{6}$ & $1.0\cdot10^{-6}$\\
Figure \ref{f:3} &$0.0125$& $0.00\dot{3}$ & $1.0\cdot10^{6}$& $3.0\cdot 10^{6}$ & $1.0\cdot10^{-6}$\\
Figure \ref{f:gs1} &$0.04$ & $0.005$ & $2.0\cdot 10^{6}$&$2.0\cdot10^{6}$ &$1.0\cdot10^{-6}$\\
Figure \ref{f:gs2} &$0.01$ & $0.01$ & $2.0\cdot 10^{6}$&$2.0\cdot10^{6}$ & $1.0\cdot10^{-6}$\\
Figure \ref{f:gs3} &$0.01$ & $0.01$ & $2.0\cdot 10^{6}$& $2.0\cdot10^{6}$ &$1.0\cdot10^{-6}$\\
\hline
\end{tabular}
\end{center}
\caption{Parameter values for simulations.}
\label{t:param}
\end{table}

\subsubsection{Binary images}
\begin{figure}[!h]
\begin{center}
\subfigure{\includegraphics[width=.24\textwidth,angle=-90]{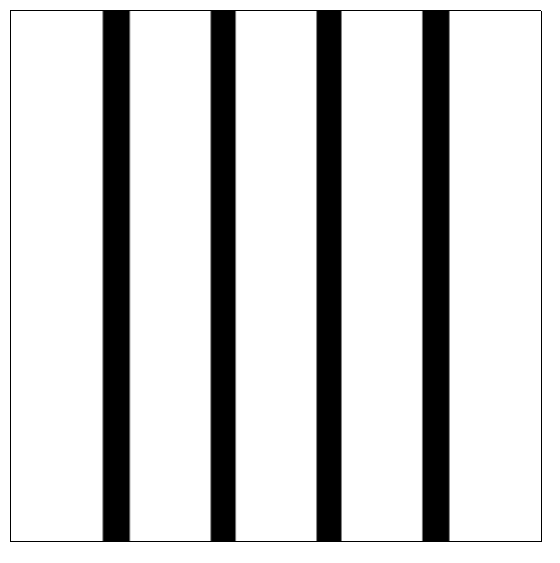}}\hspace{0mm}
\subfigure{\includegraphics[width=.24\textwidth,angle=-90]{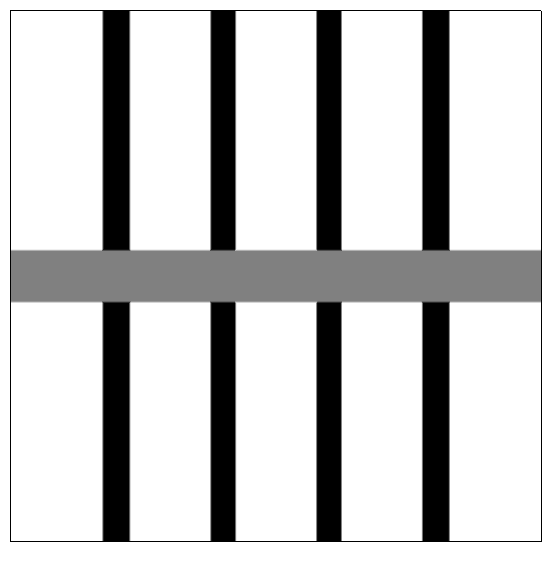}}
\hspace{0mm}
\subfigure{\includegraphics[width=.24\textwidth,angle=-90]{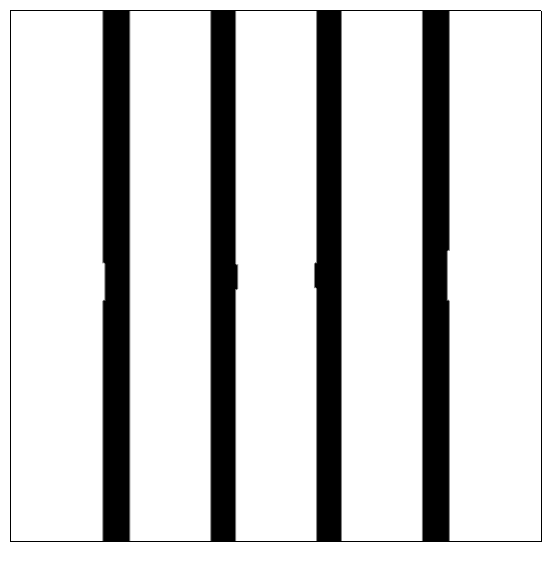}}\hspace{0mm}
\subfigure{\includegraphics[width=.24\textwidth,angle=-90]{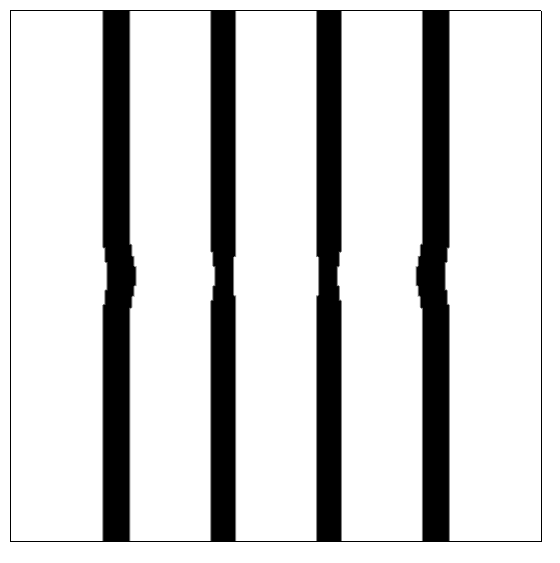}}\hspace{0mm}
\caption{Original image (first plot), damaged image (second plot),  $P(\tilde{u}_h)$ obtained using the double obstacle (third plot), $P(\tilde{u}_h)$ obtained using the quartic double well potential (fourth plot).}
\label{f:comp}
\end{center}
\end{figure}

In Figure \ref{f:comp} we compare the solution obtained using the double obstacle potential (third plot) with the solution obtained using the classical smooth quartic double well function $W_{\mathrm{qu}}(s) = (s^{2}-1)^{2}$ (fourth plot). In the first and second plots we display the original and damaged images respectively. In this example we see that the double obstacle potential outperforms the double well potential in its reconstruction of the damaged image.

\begin{figure}[!h]
\begin{center}
\subfigure{\includegraphics[width=.3\textwidth,angle=-90]{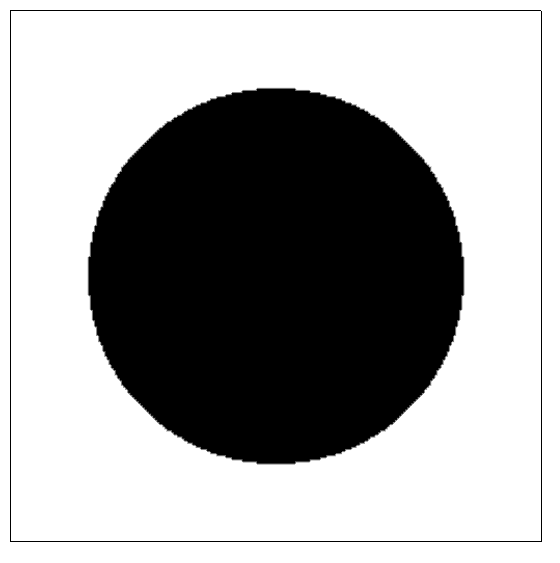}}\hspace{0mm}
\subfigure{\includegraphics[width=.3\textwidth,angle=-90]{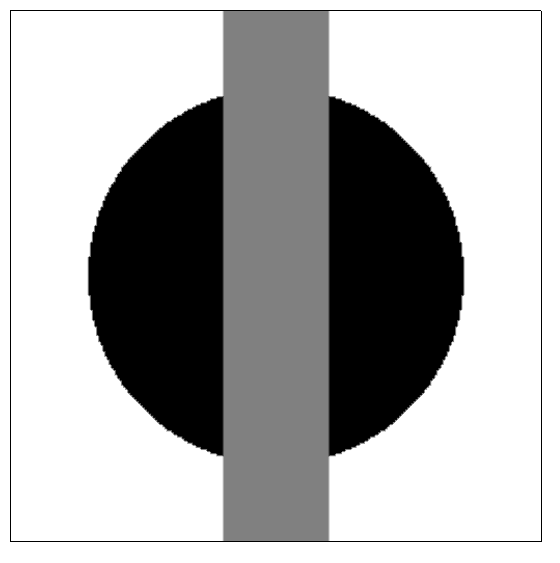}}\hspace{0mm}
\subfigure{\includegraphics[width=.3\textwidth,angle=-90]{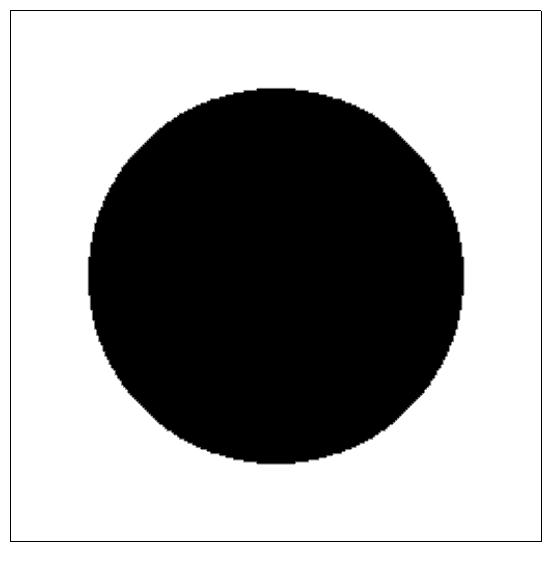}}\hspace{0mm}
\caption{Original image (left plot), damaged image (centre plot), projected solution $P(\tilde{u}_h)$ (right plot).}
\label{f:1}
\end{center}
\end{figure}

\begin{figure}[!h]
\begin{center}
\subfigure{\includegraphics[width=.3\textwidth,angle=-90]{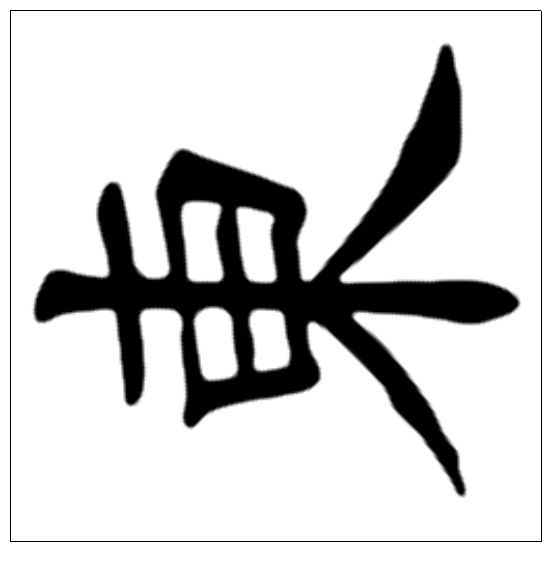}}\hspace{0mm}
\subfigure{\includegraphics[width=.3\textwidth,angle=-90]{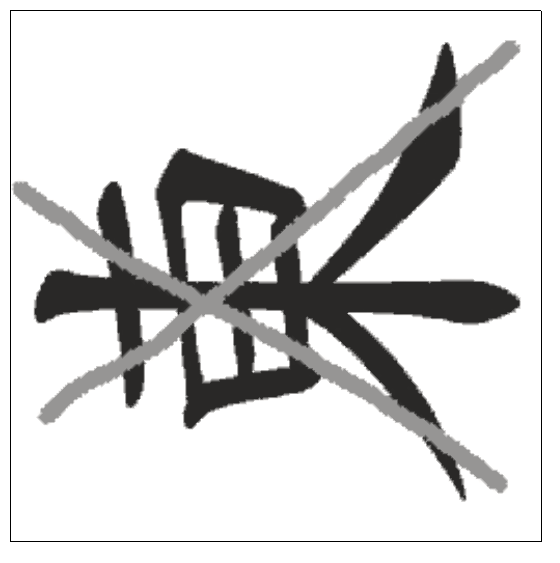}}\hspace{0mm}
\subfigure{\includegraphics[width=.3\textwidth,angle=-90]{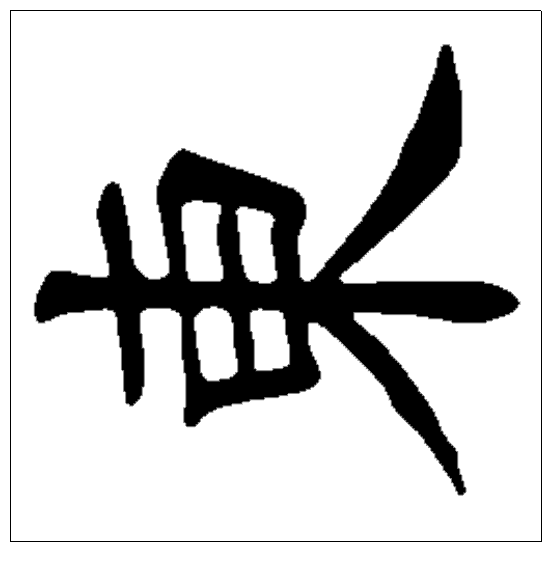}}\hspace{0mm}
\caption{Original image (left plot), damaged image (centre plot), projected solution $P(\tilde{u}_h)$ (right plot).}
\label{f:2}
\end{center}
\end{figure}

\begin{figure}[!h]
\begin{center}
\subfigure{\includegraphics[width=.3\textwidth,angle=-90]{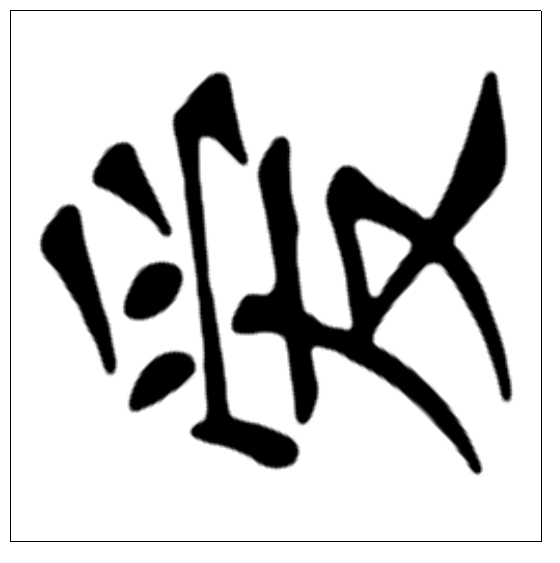}}\hspace{0mm}
\subfigure{\includegraphics[width=.3\textwidth,angle=-90]{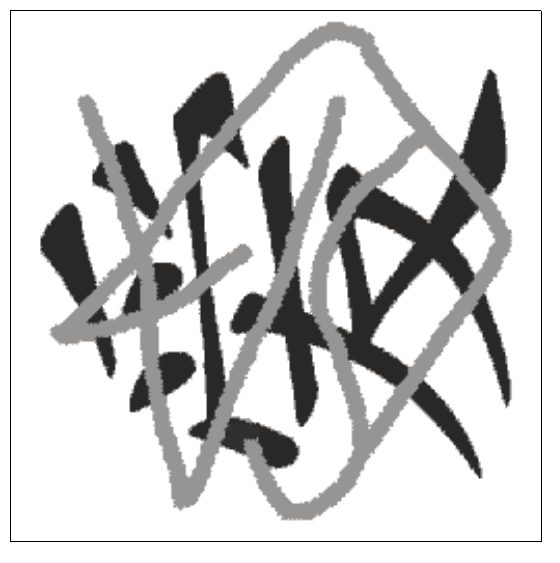}}
\subfigure{\includegraphics[width=.3\textwidth,angle=-90]{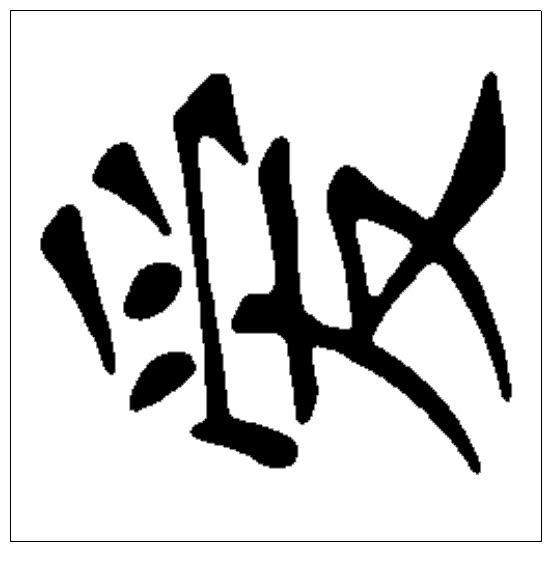}}\hspace{0mm}
\caption{Original image (left plot), damaged image (centre plot), projected solution $P(\tilde{u}_h)$ (right plot).}
\label{f:3}
\end{center}
\end{figure}

In Figures~\ref{f:1} - \ref{f:3} we display results for three binary images; in each figure we display the original image in the plot on the left, the damaged image in the centre plot, and the projected image $P(\tilde{u}_{h})$ in the plot on the right.

\begin{figure}[!h]
\begin{center}
\subfigure{\includegraphics[width=.3\textwidth,angle=0]{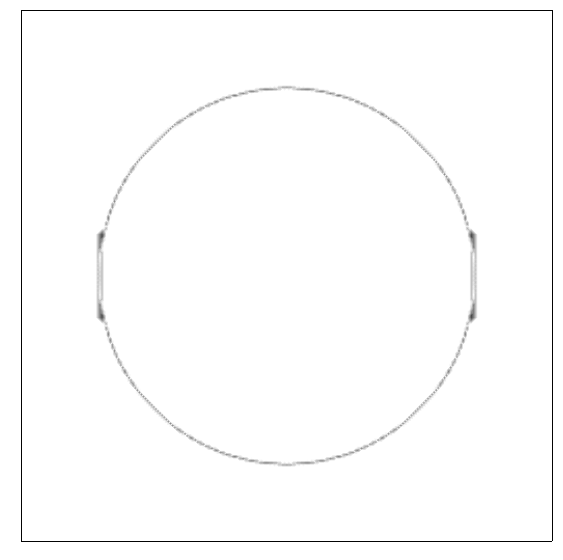}}\hspace{0mm}
\subfigure{\includegraphics[width=.3\textwidth,angle=0]{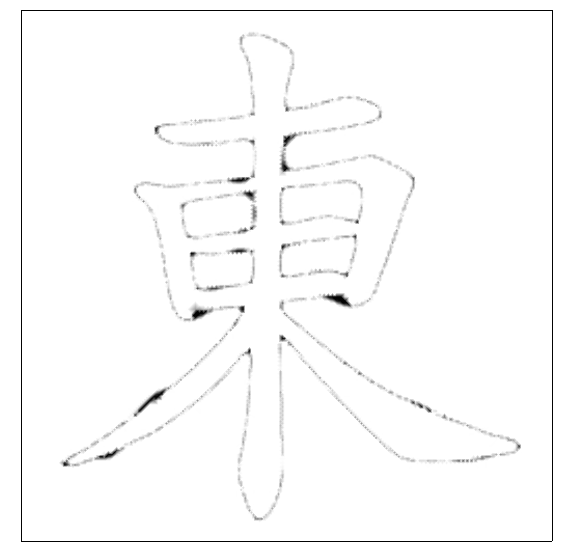}}
\subfigure{\includegraphics[width=.3\textwidth,angle=0]{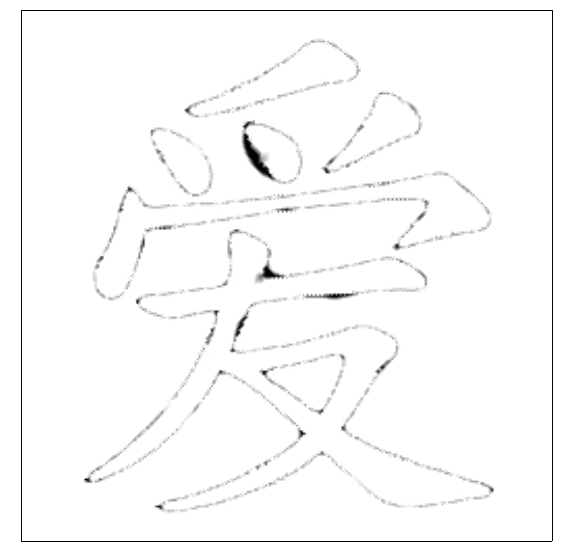}}\hspace{0mm}
\subfigure{\includegraphics[width=.07\textwidth,angle=0]{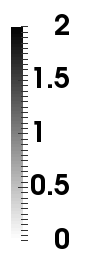}}\hspace{0mm}
\caption{The reconstruction error $|I-\tilde{u}_h|$ between the original undamaged image $I$ and the computed recovery image $\tilde{u}_h$.}
\label{f:4}
\end{center}
\end{figure}

In Figure~\ref{f:4}, we plot the error $\abs{I - \tilde{u}_{h}}$, where $I$ denotes the original undamaged image, for the three images in Figures~\ref{f:1} - \ref{f:3}.

\subsection{Grayscale images}
\begin{figure}[h]
\centering
\subfigure{\includegraphics[width = 0.32\textwidth, angle = 0]{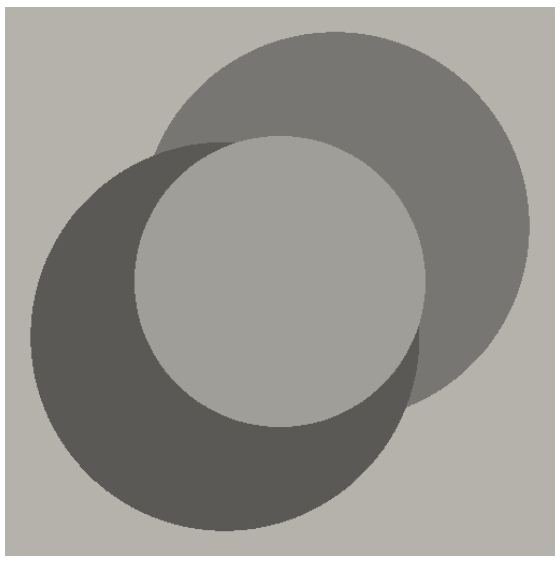}} 
\subfigure{\includegraphics[width = 0.32\textwidth, angle = 0]{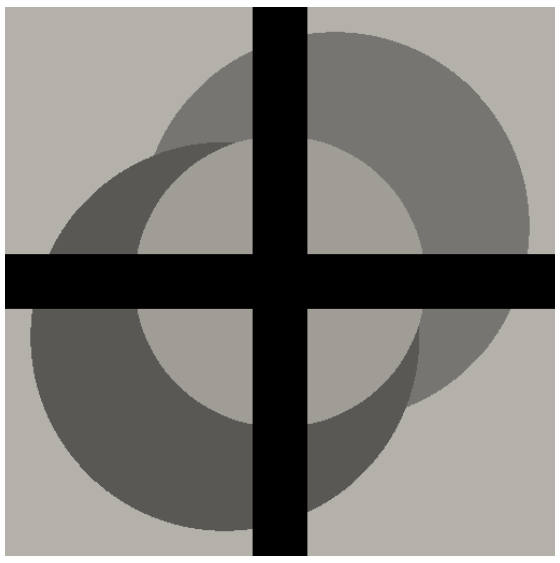}} 
\subfigure{\includegraphics[width = 0.32\textwidth, angle = 0]{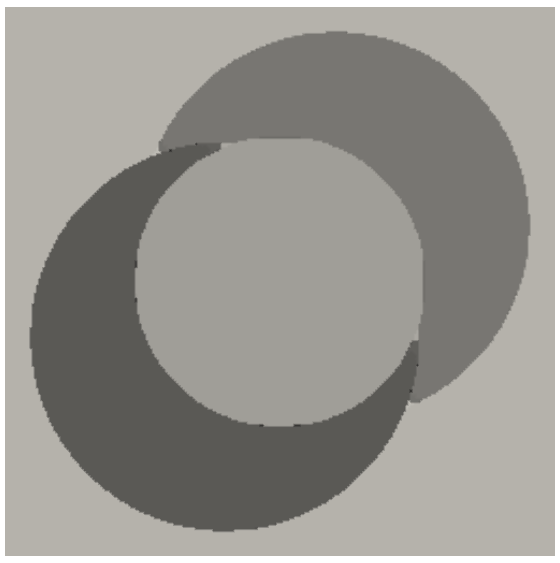}} 
\caption{original image (left), damaged image (centre), reconstructed image (right).}
\label{f:gs1}
\end{figure}

\begin{figure}[h]
\centering
\subfigure{\includegraphics[width = 0.32\textwidth, angle = 0]{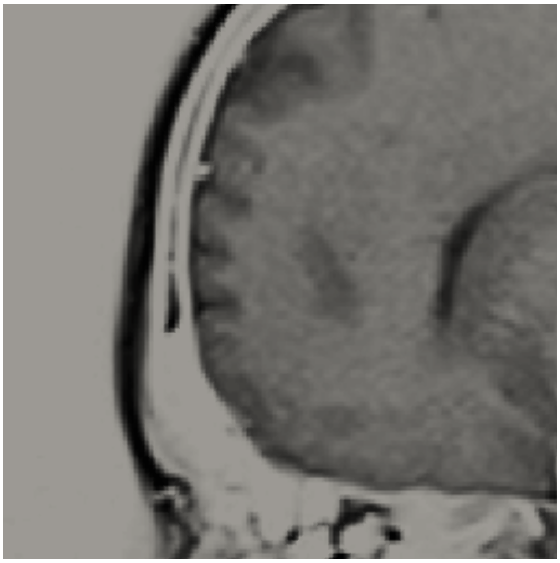}} 
\subfigure{\includegraphics[width = 0.32\textwidth, angle = 0]{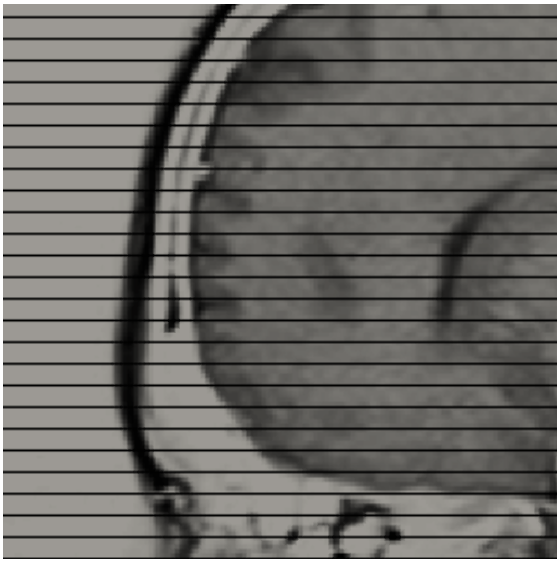}} 
\subfigure{\includegraphics[width = 0.32\textwidth, angle = 0]{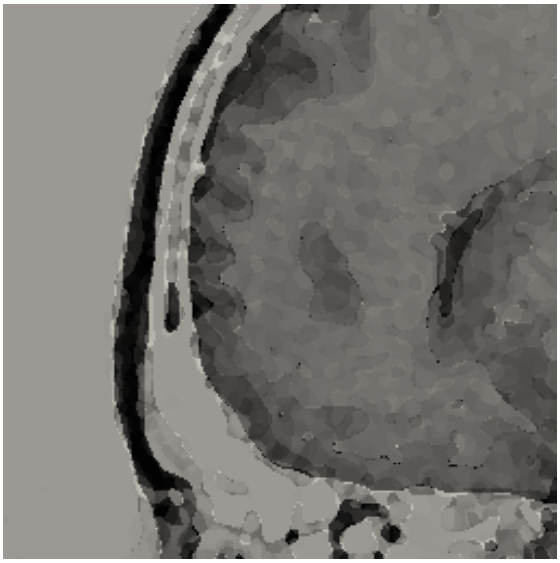}} \\
\subfigure{\includegraphics[width = 0.32\textwidth, angle = 0]{brain_id_w.png}} 
\subfigure{\includegraphics[width = 0.32\textwidth, angle = 0]{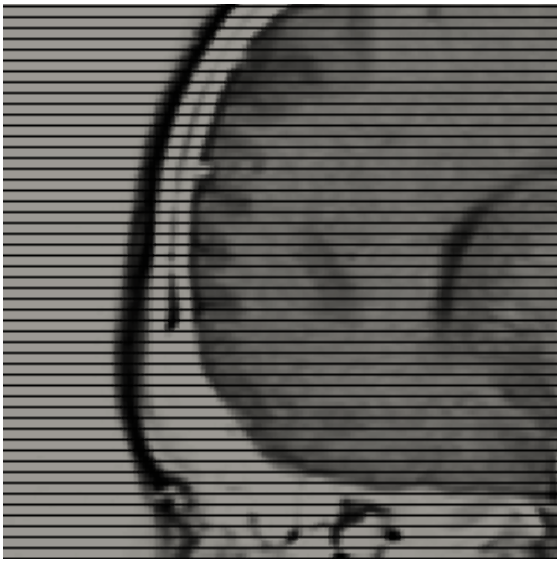}} 
\subfigure{\includegraphics[width = 0.32\textwidth, angle = 0]{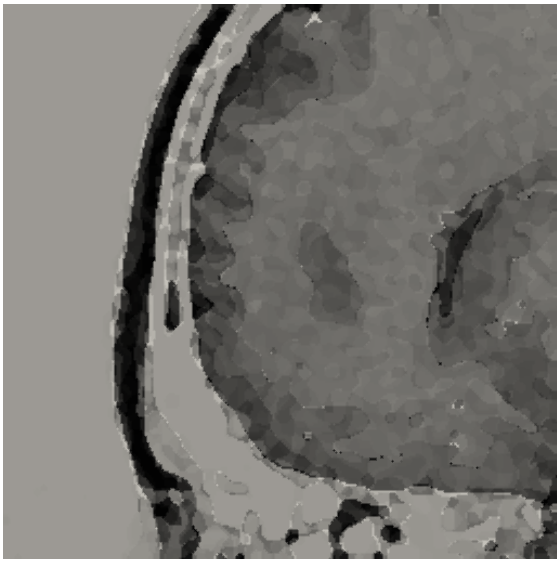}} \\
\subfigure{\includegraphics[width = 0.32\textwidth, angle = 0]{brain_id_w.png}} 
\subfigure{\includegraphics[width = 0.32\textwidth, angle = 0]{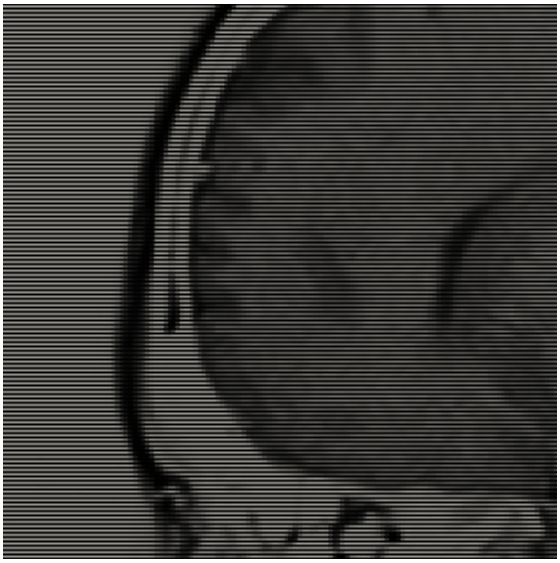}} 
\subfigure{\includegraphics[width = 0.32\textwidth, angle = 0]{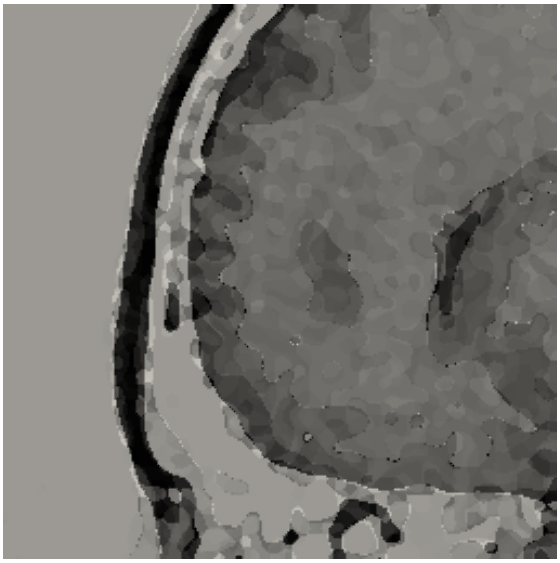}} 
\caption{original image (left), damaged image (centre) reconstructed image (right) with $10\%$ (top row), $20\%$ (middle row) and $50\%$ (bottom row) damage.}
\label{f:gs2}
\end{figure}

\begin{figure}[h]
\centering
\subfigure{\includegraphics[width = 0.32\textwidth, angle = 0]{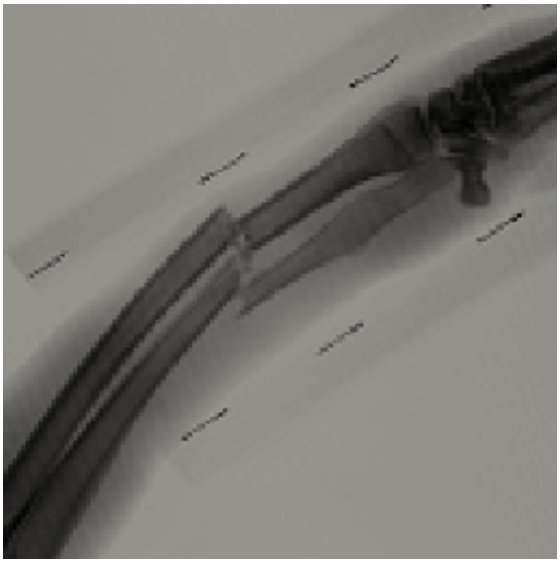}} 
\subfigure{\includegraphics[width = 0.32\textwidth, angle = 0]{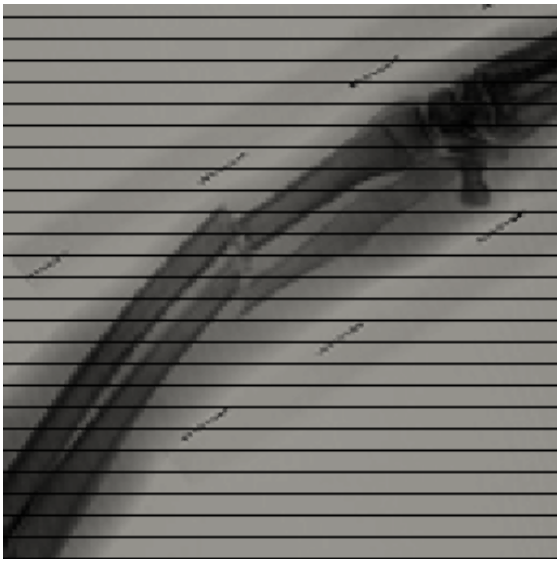}} 
\subfigure{\includegraphics[width = 0.32\textwidth, angle = 0]{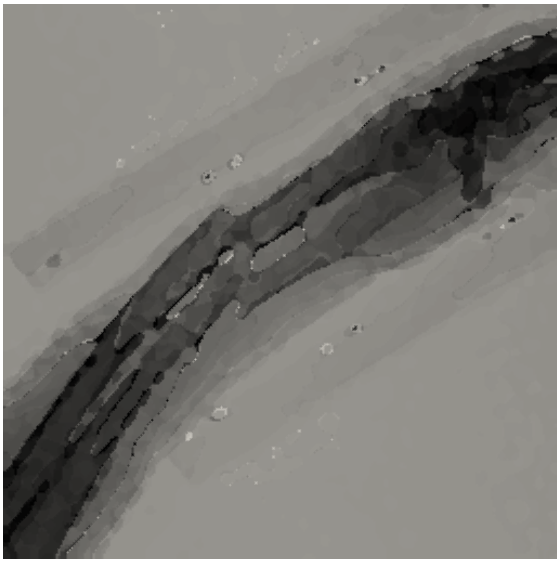}} \\
\subfigure{\includegraphics[width = 0.32\textwidth, angle = 0]{bone_id_w.png}} 
\subfigure{\includegraphics[width = 0.32\textwidth, angle = 0]{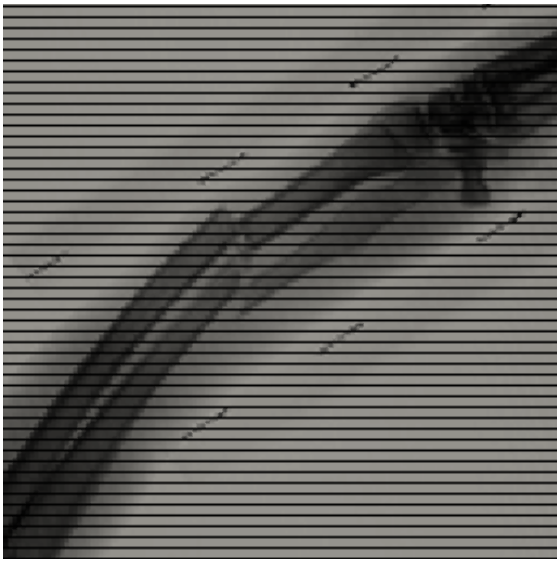}} 
\subfigure{\includegraphics[width = 0.32\textwidth, angle = 0]{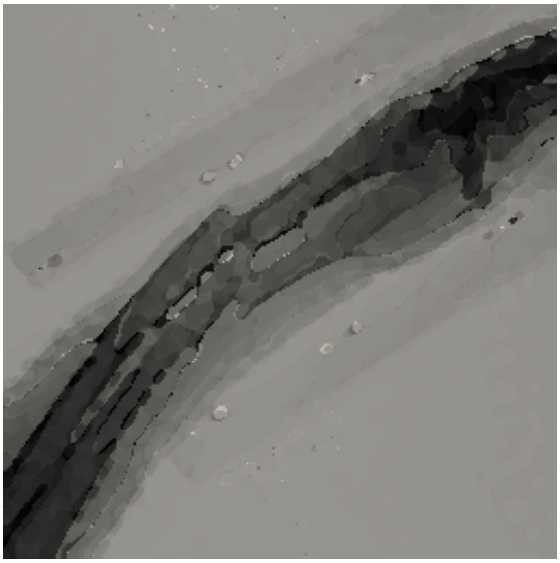}} \\
\subfigure{\includegraphics[width = 0.32\textwidth, angle = 0]{bone_id_w.png}} 
\subfigure{\includegraphics[width = 0.32\textwidth, angle = 0]{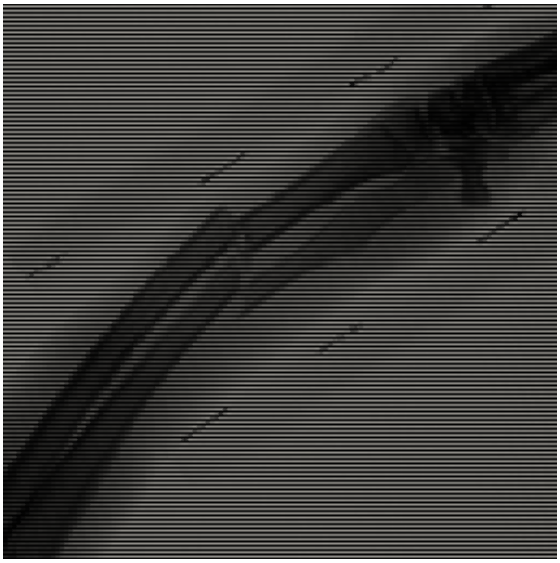}} 
\subfigure{\includegraphics[width = 0.32\textwidth, angle = 0]{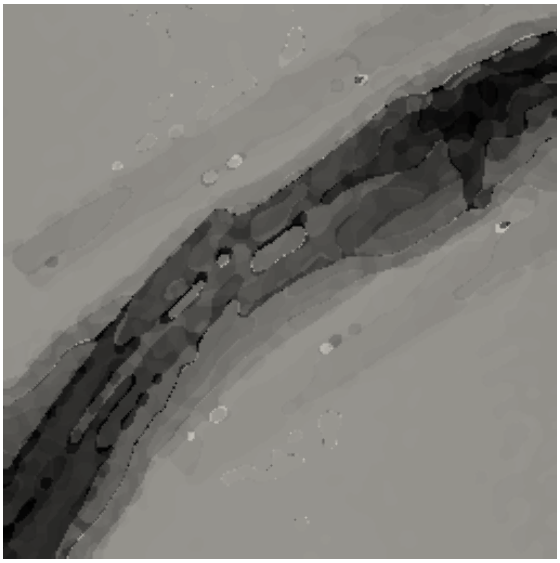}} 
\caption{original image (left), damaged image (centre), reconstructed image (right) with $10\%$ (top row), $20\%$ (middle row) and $50\%$ (bottom row) damage.}
\label{f:gs3}
\end{figure}

In \cite{sb} the authors consider a generalization of Cahn--Hilliard inpainting for grayscale images, where the grayscale image $I$ is split bitwise into $K$ channels
\begin{align*}
I(x) \approx \sum_{k=1}^K I_k 2^{-(k-1)},
\end{align*}
and the Cahn--Hilliard inpainting approach is applied to each binary channel $I_k$.  Then, at the end of the process, the channels are assembled.  In Figures \ref{f:gs1} - \ref{f:gs3} we apply this technique, with $K=8$, to some grayscale images.  We consider a simple grayscale image in Figure \ref{f:gs1} consisting of four bulk regions.  The original image (left plot), the damaged image (centre plot), and the projected solution $P(\tilde{u}_h)$ (right plot) are displayed.  Here $P(\tilde{u}_h)$ is obtained by projecting the solution from each channel onto the set of binary images and then assembling the resulting projections. 
For these results we replaced the tolerance of $5.0\cdot 10^{-6}$ in the stopping criteria in (\ref{eq:sc}) with $1.0\cdot 10^{-7}$, and the adaptivity criteria was amended so that a fine mesh is constructed at the boundaries of the bulk regions. 

We conclude this paper with the reconstructions of grayscale medical images: Figure \ref{f:gs2} is a portion of a human patient's brain MRI scan, and Figure \ref{f:gs3} is a veterinarian X-ray of a dog's broken leg.  These images are provided free for reuse and modification by www.pixabay.com.  

Each figure contains the original image (left plot), the damaged image (centre plot) and the projected reconstructed image $P(\tilde{u}_h)$ (right plot) with $10\%$ (top row), $20\%$ (middle row), and $50\%$ (bottom row) damage.  For these results we used a uniform mesh with $h_{max}=5.5\cdot 10^{-3}$.

\section*{Acknowledgements}
V. Styles acknowledges the support of the Regensburger Universit\"{a}tsstiftung Hans Vielberth.  The authors would like to thank the anonymous referees for their careful reading and suggestions which have improved the quality of the manuscript.

\bibliographystyle{plain}
\bibliography{CHpaint}

\end{document}